\documentclass[preprint,sort&compress,12pt]{elsarticle}
\usepackage{lipsum}
\makeatletter
\def\ps@pprintTitle{%
 \let\@oddhead\@empty
 \let\@evenhead\@empty
 \def\@oddfoot{}%
 \let\@evenfoot\@oddfoot}
\makeatother

\usepackage[utf8]{inputenc}
\usepackage{algpseudocode}
\usepackage{algorithm}
\usepackage[margin=1in]{geometry}
\usepackage{bold-extra}
\usepackage{xcolor}
\usepackage{adjustbox}
\usepackage{appendix}
\usepackage{lineno}

\usepackage{amssymb}
\usepackage{amsthm}
\usepackage{amsmath}
\usepackage{graphicx}
\usepackage{verbatim}
\usepackage{times}
\usepackage[titles]{tocloft}
\usepackage{enumerate}
\usepackage{bm}
\usepackage{sidecap}
\usepackage{epstopdf}
\usepackage{pdfpages}
\usepackage[colorinlistoftodos]{todonotes}
\usepackage{soul}
\usepackage{multicol}
\usepackage{float}
\usepackage{multibib}
\usepackage{geometry}
\usepackage[small,compact]{titlesec}
\usepackage[normal,bf,labelsep=colon,skip=2pt]{caption}
\usepackage{setspace}
\usepackage{subfig}
\usepackage{tocloft}
\usepackage{url}
\usepackage{floatflt}
\usepackage{array}
\usepackage{pgfgantt}
\usepackage{sidecap}
\usepackage{cite}
\usepackage{hyperref}

\numberwithin{equation}{section}
\numberwithin{figure}{section}
\numberwithin{table}{section}

\graphicspath{{./plots/}}

\DeclareMathOperator*{\argmin}{arg\,min}

\newcommand{\bR}{\mathbb R}
\newcommand{\fb}{\mathbf f}
\newcommand{\xb}{\mathbf x}

\newcommand*\diff{\mathop{}\!\mathrm{d}}

\usepackage[colorinlistoftodos]{todonotes}

\begin{document}

\begin{frontmatter}
\title{CEBoosting: Online Sparse Identification of Dynamical Systems with Regime Switching by Causation Entropy Boosting}

\author[1]{Chuanqi Chen}\ead{cchen656@wisc.edu}
\author[2]{Nan Chen}\ead{chennan@math.wisc.edu}
\author[1]{Jin-Long Wu}\ead{jinlong.wu@wisc.edu}
\address[1]{Department of Mechanical Engineering, University of Wisconsin–Madison, Madison, WI 53706}
\address[2]{Department of Mathematics, University of Wisconsin–Madison, Madison, WI 53706}

\begin{abstract}
Regime switching is ubiquitous in many complex dynamical systems with multiscale features, chaotic behavior, and extreme events. In this paper, a causation entropy boosting (CEBoosting) strategy is developed to facilitate the detection of regime switching and the discovery of the dynamics associated with the new regime via online model identification. The causation entropy, which can be efficiently calculated, provides a logic value of each candidate function in a pre-determined library. The reversal of one or a few such causation entropy indicators associated with the model calibrated for the current regime implies the detection of regime switching. Despite the short length of each batch formed by the sequential data, the accumulated value of causation entropy corresponding to a sequence of data batches leads to a robust indicator. With the detected rectification of the model structure, the subsequent parameter estimation becomes a quadratic optimization problem, which is solved using closed analytic formulae. Using the Lorenz 96 model, it is shown that the causation entropy indicator can be efficiently calculated, and the method applies to moderately large dimensional systems. The CEBoosting algorithm is also adaptive to the situation with partial observations. It is shown via a stochastic parameterized model that the CEBoosting strategy can be combined with data assimilation to identify regime switching triggered by the unobserved latent processes. In addition, the CEBoosting method is applied to a nonlinear paradigm model for topographic mean flow interaction, demonstrating the online detection of regime switching in the presence of strong intermittency and extreme events.


\end{abstract}
\end{frontmatter}

\section{Introduction}

Regime switching is ubiquitous in many complex dynamical systems in geoscience, engineering, neural science, and material science \citep{majda2016introduction, majda2006nonlinear, strogatz2018nonlinear, baleanu2011fractional, deisboeck2007complex, stelling2001book, sheard2009principles}. The switching is usually associated with sudden changes in internal states or appears when specific external forcing is exerted. Dynamical systems often display distinct behavior with regime switching. One example is the atmospheric jets, which meander in different directions when the atmosphere alternates between blocked and unblocked regimes \citep{amundson2010ice, majda2013elementary}. Similarly, an excitable medium is susceptible to finite perturbations, which triggers regime switching from a quiescent state to one with various wave patterns \citep{holden2013nonlinear, lindner2004effects, chen2019spatial}. Regime switching can also induce an increased occurrence of extreme events, leading to, for example, extreme weather and climate patterns \citep{lau2011intraseasonal, clarke2008introduction}, bursting neurons \citep{altmann2005recurrence}, or extreme ductile damages \citep{bronkhorst2021local}. Detecting regime switching and the corresponding underlying dynamics, which relies on appropriate model identification methods, has significant social and scientific impacts. Challenges in detecting regime switching are associated with the intrinsic properties of many complex dynamical systems, including high dimensionality, partial or incomplete observations, and the intermittent occurrence of rare and extreme events \citep{majda2018model, wiggins2003introduction, kalnay2003atmospheric, law2015data, branicki2013non}.

Efficient model identification has received significant attention. Both physical knowledge and observational data facilitate learning the underlying model dynamics. The model structures are often established utilizing physical intuitions for traditional knowledge-based model identification. The primary process then becomes the estimation of model parameters. Linear models are natural candidates for simple problems \citep{freedman2009statistical, yan2009linear} and can potentially be skillful for short-term forecasts. Other families of models with pre-determined structures, such as the physics-constrained nonlinear regression models \citep{majda2012physics, harlim2014ensemble} and conditional Gaussian nonlinear systems \citep{chen2018conditional, chen2022conditional}, are alternative nonlinear models aiming to capture specific underlying dynamical features. On the other hand, recent progress has been made in data-driven model identification. Data-driven reduced-order models have been widely used in scientific and engineering applications \citep{ahmed2021closures, hijazi2020data, lin2021data, peherstorfer2015dynamic}. Sparse model discovery methods also appear as advanced model identification tools that allow automatic learning of the model structure and parameters from data and lead to nonlinear models with parsimonious structures via sparse regression  \citep{brunton2016discovering, rudy2017data, fasel2022ensemble, schaeffer2013sparse, billings2007sparse, mojgani2022discovery, quade2018sparse, chen2020learning}. With a limited amount of indirect data, derivative-free optimization methods~\citep{xiao2016quantifying, zhang2022ensemble, schneider2022ensemble} have been explored as model identification tools. In addition, non-parametric and machine learning models have been built to characterize complex dynamical systems \citep{pawar2020data, moosavi2015efficient, san2018extreme, beck2021perspective, wang2017physics, wu2018physics}.

Among various model discovery approaches, online model identification is a particularly useful method in practice, which sequentially determines model structure and estimates model parameters when new observation arrives \citep{rong2006sequential, lombaerts2009nonlinear, kopsinis2010online, kalouptsidis2011adaptive,  chen2014system, srinivasan2019sequential}. It is the primary model identification strategy in many geophysical and engineering problems, where the limited amount of historical data is insufficient to robustly discover the underlying dynamics. It should be noted that online model identification can be further combined with data assimilation to handle noisy observations or recover the unobserved state variables in the situation with partial observations \citep{gottwald2021supervised, wikner2021using, schneider2021learning}. However, unlike the online parameter estimation that can be efficiently addressed by standard filtering methods~\citep{kalman1960new}, the lack of knowledge about the proper model structure poses a unique challenge in online model identification. The sequentially arriving data plays a vital role in progressively rectifying the model and reducing the uncertainty in the identified system. Although existing system identification methods can identify the proper model structure via promoting sparsity in the offline setting with fitting a model to abundant data, promoting sparsity relies on the model fitting may eliminate some important model structures in the context of sequential learning. To address such a challenge, we incorporate causation entropy to achieve robust online model identification. As regime switching often occurs and completes within a short transient period, developing suitable online identification methods for discovering regime switching exploiting transition data is essential with practical importance.

In this paper, a causation entropy boosting (CEBoosting) strategy is developed. It is incorporated into an online model identification method to detect regime switching and discover the nonlinear dynamics associated with the new regime. Different from many existing sparse model identification algorithms, such as those relying on LASSO (least absolute shrinkage and selection operator) regression~\citep{tibshirani1996regression,schneider2022ensemble} or thresholding~\citep{brunton2016discovering,schaeffer2013sparse}, the method developed here separates the estimation of model parameters from the recurrent identification of nonlinear model structure. Such a separation allows using closed analytic formulae for the entire online learning algorithm, and therefore, the overall computational cost is significantly reduced. In this new strategy, causation entropy \citep{elinger2020information, elinger2021causation} is utilized to provide a logic value (i.e., true or false) of each candidate function in a pre-determined library throughout the online learning process. By examining the causation entropy on the newly arrived data, the reversal of one or a few such causation entropy indicators associated with the model calibrated for the current regime implies the detection of regime switching. In other words, the causation entropy indicator, which can be efficiently calculated, is employed to decide if the existing terms in the current model need to be rectified and if the system demands additional terms as a response to regime switching. Note that the sequential data in online learning is collected within a short time window to form a batch of time series, which is utilized to compute the causation entropy. As each batch contains a short amount of data, it may embody only part of the dynamical properties. Nevertheless, as time evolves, the accumulated value of causation entropy corresponding to a sequence of batches leads to a robust indicator of the model structure in response to regime switching. The concept of accumulating causation entropy calculated from sequential data relates to the statistical method of bagging~\citep{breiman1996bagging}. With the detected rectification of the model structure, the subsequent parameter estimation becomes a quadratic optimization problem, which is solved using closed analytic formulae. For multiple times of regime switching, a summation of residual models is calibrated, which relates to the statistical method of boosting~\citep{freund1996experiments, friedman2002stochastic, schapire2013boosting, chen2015xgboost}.

The proposed new strategy has several unique features. First, causation entropy takes into account the interdependence between all the candidate functions in the pre-determined library, and therefore, it can eliminate the superficial causal relationship. Model identification exploiting the causation entropy has been shown to reach a higher selection accuracy than LASSO regression or elastic net \citep{elinger2021causation}. The causation-based learning approach also indicates robust results in the presence of indirect coupling between features and stochastic noise \citep{quinn2015directed}, which are crucial features of complex systems. Second, causation entropy is only utilized to indicate the terms that need to be added or removed from the existing model. In other words, although computing the exact value of the causation entropy is challenging, closed analytic formulae are available for efficiently approximating this causation entropy indicator, which allows an effective detection of the model structure. Third, the parameter estimation only needs to be carried out after the model structure is entirely determined. Therefore, the overall computational cost is reduced compared with applying LASSO regression, which requires detecting the model structure and estimating model parameters simultaneously for each batch of data. It is worth highlighting that the causation entropy indicator is easy to calculate and applicable to moderately large dimensional systems. The method developed here is also adaptive to the situation with partial observations, where utilizing data assimilation to recover the unobserved state variables can be incorporated into the learning process for identifying regime switching resulting from the latent processes. In addition, the method is not limited to the Gaussian data. It can be applied to dynamical regimes with strong intermittency and extreme events. Applications to nonlinear dynamical systems with moderately large dimensions, partial observations, and extreme events are all studied in the paper.

The remainder of the paper is organized as follows. The development of the online identification method utilizing the CEboosting strategy is presented in Section~\ref{sec:methodology}. Section~\ref{sec:results} includes four test cases. In addition to a standard chaotic model as a proof-of-concept, the other three test cases emphasize the method applying to systems with moderately large dimensions, partial observations, and extreme events, respectively. The paper is concluded in Section~\ref{sec:conclusion}.

\section{Methodology}
\label{sec:methodology}

The online sparse identification method aims to (i) detect regime switching of dynamical systems via causation entropy and (ii) determine the resulting dynamics after the regime switching. The dynamical system has the following general form:
\begin{equation}
\label{eq:system}
    \dot{\textbf{x}}(t) = \textbf{f}(\textbf{x}(t)) + \boldsymbol{\sigma} \dot{\textbf{W}}(t),
\end{equation}
where $\textbf{x}(t)=[x_1(t),x_2(t),\ldots,x_p(t)]^\top \in \bR^p$ is the multi-dimensional state variable and $\dot{\textbf{x}}(t) \in \bR^p$ is the associated temporal derivative, $\textbf{f}:\ \bR^p \mapsto \bR^p$ is a vector-valued nonlinear function (i.e., vector field) of the state variable, $\dot{\textbf{W}}(t) \in \bR^q$ is a white noise vector and $\boldsymbol{\sigma} \in \bR^{p \times q}$ is a matrix of noise magnitudes. In the absence of random noise forcing, $\boldsymbol{\sigma}$ becomes a zero matrix. Assume that the regime switching occurs at $t=t_s$, i.e., the vector field $\textbf{f}$ of the original dynamical system changes to $\textbf{f}^*$. The goal of this work is to detect such a regime switching and to identify the new model after the regime switching.

As regime switching often results from a sudden change of a small number of the model parameters or specific components of the model structure, the residual model $\delta \textbf{f}=\textbf{f}^*-\textbf{f}$ typically has a sparse structure that can be calibrated using relatively short data, which is precisely the case of the online identification problems. Therefore, instead of learning the entire model associated with the new regime, the focus is to estimate the residual part $\delta \textbf{f}=\textbf{f}^*-\textbf{f}$. Once the residual part $\delta \textbf{f}$ is identified, it is then added to the existing model that provides the new system as a response to the regime switching.

The limited amount of data is assumed to arrive sequentially in the form of batches. The $k$-th batch represents data $\textbf{x}(t)$ (or a subset of the vector $\textbf{x}(t)$ in the partial observation case) for $t \in [t_{B_k},t_{B_{k+1}})$. It should be noted that, as $t_s$ is typically unknown in practice, the identification algorithm usually does not start from $t_s$ (namely the left point of the first interval $t_{B_1}\neq t_s$). Yet, for the simplicity of presentation, $t_{B_1}$ is chosen to be $t_s$ for the numerical examples in this work. This will not affect the identification algorithm as applying the algorithm to those batches prior to $t_s$ will not indicate regime switching. But this setup facilitates counting for the length of the data that is needed to detect the regime switching once it occurs at $t_s$.

Denote by $\boldsymbol{\Phi}=[\phi_1,\phi_2,...,\phi_N]^\top$ a vector containing all candidate basis functions, which are knowledge-based and are pre-determined. Each $\phi_n$ in $\boldsymbol{\Phi}$ is a scalar-valued function $\phi_n:=\phi_n(\textbf{x})$ that gives a map $\bR^p \mapsto \bR$. The representation of $\delta \textbf{f}=\textbf{f}^*-\textbf{f}$ is approximated by a linear combination of these basis functions:
\begin{equation}
\label{eq:delta_f}
    \delta f_i = \sum_{n=1}^N \xi_{in} \phi_n,
\end{equation}
where $\delta f_i$ is the $i$-th scalar component of $\delta \textbf{f}$. A sparse representation of \eqref{eq:delta_f} means that most of the coefficients $\xi_{in}$ are zeros in the identified model. To obtain such a sparse representation of \eqref{eq:delta_f}, a CEBoosting method is developed to effectively determine which basis functions should take non-zero coefficients.

\subsection{Causation Entropy}
Causation entropy is based on the general concept of conditional mutual information~\citep{wyner1978definition}. Other popular information measures that also build on conditional mutual information include directed information~\citep{massey1990causality} and transfer entropy~\citep{schreiber2000measuring}, which has found applications in many areas, e.g., turbulence modeling~\citep{lozano2020causality,lozano2022information} and neurosciences~\citep{vicente2011transfer}, with a comprehensive review in~\citep{bossomaier2016transfer}. The idea of utilizing the causation entropy to detect the influence between different variables has been studied in \citep{sun2014causation, branchini2022causal, sun2015causal, elinger2020information, elinger2021causation, chen2023causality}. It can be naturally applied to the context of system identification. As the white noise does not explicitly contribute to the causal relationship, the calculation of the causation entropy mainly focuses on the candidate functions that consist of the deterministic part of the dynamics, namely the functions $\textbf{f}$ in \eqref{eq:system}. To this end, consider the deterministic part of \eqref{eq:system}:
\begin{equation}
\label{eq:discretized_system}
    \begin{bmatrix}
        \dot{x}_1(t) \\
        \dot{x}_2(t) \\
        \vdots \\
        \dot{x}_p(t)
    \end{bmatrix}
     =
    \begin{bmatrix}
        \xi_{1,1} \quad \cdots \quad \xi_{1,N} \\
        \xi_{2,1} \quad \cdots \quad \xi_{2,N} \\
        \vdots \quad \ddots \quad \vdots \\
        \xi_{p,1} \quad \cdots \quad \xi_{p,N}
    \end{bmatrix}
    \begin{bmatrix}
        \phi_1(x_1(t),...,x_p(t)) \\
        \phi_2(x_1(t),...,x_p(t)) \\
        \vdots \\
        \phi_N(x_1(t),...,x_p(t)) \\
    \end{bmatrix}  = \boldsymbol{\Xi}\boldsymbol{\Phi}.
\end{equation}
The causation entropy $C_{\phi_n \rightarrow \dot{x}_i|[\boldsymbol{\Phi} \setminus \phi_n]}$ is utilized to quantify the contribution from the candidate function $\phi_n$ to the dynamics $\dot{x}_i$ (i.e., the time derivative of the $i$-th state variable: $x_i$) conditioned on the remaining candidate functions $\boldsymbol{\Phi} \setminus \phi_n$, namely all the candidate functions except $\phi_n$. This causation entropy reflects the causal influence of $\phi_n$ to the dynamics $\dot{x}_i$, and we enforce $\xi_{in}=0$ if the causation entropy is small. Repeating this procedure over all $n=1,\ldots, N$ and $i=1,\ldots,p$ to form the matrix $\boldsymbol{\Xi}$. As only a few candidate functions will have the actual causal influence on the dynamics, the matrix $\boldsymbol{\Xi}$ is expected to have a sparse structure. The causation entropy $C_{\phi_n \rightarrow \dot{x}_i|[\boldsymbol{\Phi} \setminus \phi_n]}$ is defined as follow:
\begin{equation}
\label{eq:C}
    C_{\phi_n \rightarrow \dot{x}_i|[\boldsymbol{\Phi} \setminus \phi_n]} = H(\dot{x}_i | [\boldsymbol{\Phi} \setminus \phi_n]) - H(\dot{x}_i | \boldsymbol{\Phi}),
\end{equation}
where $H(\cdot|\cdot)$ is the conditional entropy, which is defined as:
\begin{equation}
\label{eq:cond_entropy}
    H(V|U) = \int_u \int_v p(u,v)\log(p(v|u))\diff v\diff u,
\end{equation}
where $p(u,v)$ is the corresponding probability density function (PDF) that can be determined by a histogram from the time series assuming ergodicity. On the right-hand side of \eqref{eq:C}, the difference between the two conditional entropies indicates the information in $\dot{x}_i$ contributed by the specific function $\phi_n$ given the contributions from all the other functions in the library $\bm{\Phi}$. Thus, it tells if $\phi_n$ provides additional information to $\dot{x}_i$. It is worth highlighting that the causation entropy in \eqref{eq:C} is fundamentally different from directly computing the correlation between $\dot{x}_i$ and $\phi_n$, as the causation entropy also considers the influence of the other library functions. If both $\dot{x}_i$ and $\phi_n$ are caused by a common factor $\phi_{m}$, then $\dot{x}_i$ and $\phi_n$  can be highly correlated. Yet, in such a case, the causation entropy $C_{\phi_{n} \rightarrow \dot{x}_i \mid\left[\boldsymbol{\Phi}  \backslash {\phi}_{n}\right]}$ will be zero as $\phi_n$ is not the causation of $\dot{x}_i$.

In practice, the conditional entropy in \eqref{eq:cond_entropy} can involve expensive high-dimensional integrals, which is computationally challenging \citep{bellman1961dynamic}. Nevertheless, a Gaussian approximation of the PDFs inside the integrand can be utilized to calculate the causation entropy \citep{chen2023causality}. By approximating all the joint and marginal distributions as Gaussians, the causation entropy is calculated as follows:
\begin{equation}
\label{eq:C_gaussian}
\begin{aligned}
    C_{W \rightarrow U|V} &= H(U|V) - H(U|V,W) \\
    &= H(U,V) - H(V) - H(U,V,W) + H(V,W) \\
    &= \frac{1}{2}\ln(\det(\mathbf{R}_{UV})) - \frac{1}{2}\ln(\det(\mathbf{R}_{V})) - \frac{1}{2}\ln(\det(\mathbf{R}_{UVW})) + \frac{1}{2}\ln(\det(\mathbf{R}_{VW})),
\end{aligned}
\end{equation}
where $\mathbf{R}$ denotes the covariance matrix of the corresponding vector, e.g., $\mathbf{R}_{UVW}$ corresponds to the covariance matrix of the vector $[U,V,W]^\top$. The explicit expression in \eqref{eq:C_gaussian} based on the Gaussian approximation can efficiently compute the causation entropy. It allows the computation of the causation entropy with a moderately large dimension, which is typically the case for many practical situations. It is worth noting that the Gaussian approximation may lead to certain errors in computing the causation entropy if the actual distribution is highly non-Gaussian. Nevertheless, the primary goal is not to obtain the exact value of the causation entropy. Instead, it suffices to detect if the causation entropy $C_{\phi_{n} \rightarrow \dot{x}_i \mid\left[\boldsymbol{\Phi}  \backslash {\phi}_{n}\right]}$ is nonzero (or practically above a small threshold value). In most applications, if a significant causal relationship is detected in the higher-order moments, it is very likely in the Gaussian approximation. This allows us to efficiently determine the sparse model structure. The exact values of the nonzero coefficients on the right-hand side of the identified model will be calculated via a simple least square estimation to be discussed in the following. Note that the Gaussian approximation is taken directly from the statistics associated with the nonlinear time series from the underlying nonlinear model. Therefore, the Gaussian approximation still includes the nonlinear dynamical information. It is very different from linearizing a nonlinear complex system and computing the resulting Gaussian distribution. Such a Gaussian approximation of the nonlinear time series has been widely applied to compute various information measurements and lead to reasonably accurate results \citep{majda2018model, tippett2004measuring, kleeman2011information, branicki2012quantifying}. Note that, with linear and Gaussian assumptions, causation entropy has been demonstrated as equivalent to Granger causality~\citep{barnett2009granger}, which has been a popular tool for analyzing time series data for decades~\citep{shojaie2022granger}. But the focus here is more toward identifying the nonlinear models.

Below, for notation conciseness, $\mathbf{C}_{in}$ is utilized as a short-hand notation of $C_{\phi_n \rightarrow \dot{x}_i|[\boldsymbol{\Phi} \setminus \phi_n]}$. To impose sparsity into $\boldsymbol{\Phi}$ in the practical computational scenarios, a threshold $\overline{C}$ is prescribed and $\mathbf{C}_{in}=0$ is enforced when $C_{\phi_n \rightarrow \dot{x}_i|[\boldsymbol{\Phi} \setminus \phi_n]} \leq \overline{C}$. Such a threshold value is adopted mainly to exclude the small causation entropy values due to the sampling error from using a finite time series. After applying this threshold value,  a causation entropy matrix (CEM) $\mathbf{C}$ is obtained, where its $(i,n)$-th entry is given by:
\begin{equation}
\label{eq:threshold}
    \mathbf{C}_{in} =
    \begin{cases}
        0 \quad \text{if $C_{\phi_n \rightarrow \dot{x}_i|[\boldsymbol{\Phi} \setminus \phi_n]} \leq \overline{C}$}, \\
        1 \quad \text{if $C_{\phi_n \rightarrow \dot{x}_i|[\boldsymbol{\Phi} \setminus \phi_n]} > \overline{C}$}.
    \end{cases}
\end{equation}
Based on the matrix $\mathbf{C}$, the sparsity can be further enforced into $\boldsymbol{\Xi}$ by setting $\xi_{in}=0$ when $\mathbf{C}_{in}=0$. It has been demonstrated that a correct sparse model $\boldsymbol{\Xi}\boldsymbol{\Phi}$ can be obtained in the offline-learning setting where the time series is long enough \citep{chen2023causality}. However, in the online learning setting, all the covariance matrices in \eqref{eq:C_gaussian} are only estimated from a limited amount of batch data and may not provide correct information for imposing sparsity. To address this challenge in the online learning setting, the following CEBoosting algorithm is introduced.

\subsection{Causation Entropy Boosting (CEBoosting) Algorithm}

The $k$-th batch data corresponds to a time series $\textbf{x}(t)$ for $t \in [t_{B_k},t_{B_{k+1}})$. Assume that by exploiting all the past batch data $\textbf{x}(t)$ for $t \in [0,t_{B_{1}})$, the model in the current regime is estimated as:
\begin{equation}
\label{eq:null_system}
    \begin{bmatrix}
        \dot{x}_1(t) \\
        \dot{x}_2(t) \\
        \vdots \\
        \dot{x}_p(t)
    \end{bmatrix}
     =
    \begin{bmatrix}
        \xi^{(0)}_{1,1} \quad \cdots \quad \xi^{(0)}_{1,N} \\
        \xi^{(0)}_{2,1} \quad \cdots \quad \xi^{(0)}_{2,N} \\
        \vdots \quad \ddots \quad \vdots \\
        \xi^{(0)}_{p,1} \quad \cdots \quad \xi^{(0)}_{p,N}
    \end{bmatrix}
    \begin{bmatrix}
        \phi_1(x_1(t),...,x_p(t)) \\
        \phi_2(x_1(t),...,x_p(t)) \\
        \vdots \\
        \phi_N(x_1(t),...,x_p(t)) \\
    \end{bmatrix}  = \boldsymbol{\Xi}^{(0)}\boldsymbol{\Phi}.
\end{equation}

With the $1$-st new incoming batch data $\textbf{x}(t)$ for $t \in [t_{B_1},t_{B_{2}})$, the true dynamics $\dot{\textbf{x}}$ within this period can be obtained by differentiation of $\textbf{x}(t)$. On the other hand, the predicted dynamics using the current model \eqref{eq:null_system} is given by $\dot{\tilde{\textbf{x}}} = \boldsymbol{\Xi}^{(0)}\boldsymbol{\Phi}$. Then the residual dynamics $\delta \textbf{f}$ for the $1$-st batch is given by $\textbf{r} = \dot{\textbf{x}} -\dot{\tilde{\textbf{x}}}$. If there is no regime switching for this batch, then it is expected that $\delta \textbf{f} = 0$ and the estimated $\textbf{r}$ is approximately Gaussian white noise. As an analog to \eqref{eq:delta_f}, the matrix form of the residual dynamics $\textbf{r}$ can be written as:
\begin{equation}
\label{eq:residual_system}
    \begin{bmatrix}
        r_1(t) \\
        r_2(t) \\
        \vdots \\
        r_p(t)
    \end{bmatrix}
     =
    \begin{bmatrix}
        \delta\xi_{1,1} \quad \cdots \quad \delta\xi_{1,N} \\
        \delta\xi_{2,1} \quad \cdots \quad \delta\xi_{2,N} \\
        \vdots \quad \ddots \quad \vdots \\
        \delta\xi_{p,1} \quad \cdots \quad \delta\xi_{p,N}
    \end{bmatrix}
    \begin{bmatrix}
        \phi_1(x_1(t),...,x_p(t)) \\
        \phi_2(x_1(t),...,x_p(t)) \\
        \vdots \\
        \phi_N(x_1(t),...,x_p(t)) \\
    \end{bmatrix}  = \boldsymbol{\Xi}_r\boldsymbol{\Phi}.
\end{equation}
To detect regime switching, identify sparse pattern of $\boldsymbol{\Xi}_r$ and estimate the rest of parameters of $\boldsymbol{\Xi}_r$ in \eqref{eq:residual_system}, one can estimate the causation entropy $C_{\phi_n \rightarrow r_i|[\boldsymbol{\Phi} \setminus \phi_n]}$ between basis of states $\boldsymbol{\Phi}$ and residual dynamics $\textbf{r}$ according to \eqref{eq:C} and \eqref{eq:C_gaussian} and then obtain the corresponding causation entropy matrix $\mathbf{C}$ by \eqref{eq:threshold}.

Mathematically, if $\mathbf{C} = \textbf{0}$, it implies that no variable in the basis functions $\boldsymbol{\Phi}$ is significant to the residual dynamics $\textbf{r}$, which leads to $\boldsymbol{\Xi}_r=\textbf{0}$ in \eqref{eq:residual_system}. In other words, the current model $\boldsymbol{\Xi}^{(0)}\boldsymbol{\Phi}$ fits the dynamics of $1$-st batch and make residual dynamics $\textbf{r}$ behave like white noise. This means there is no regime switching ($\delta \textbf{f} =\textbf{0} $) if $\mathbf{C} = \textbf{0}$ and the model is unchanged. Otherwise, if $\mathbf{C}$ contains some non-zero terms, then a regime switching occurs ($\delta \textbf{f} \neq \textbf{0} $), and the goal is to identify the sparse structure of the $\boldsymbol{\Xi}_r$ in \eqref{eq:residual_system}. However, despite the mathematical justification, with a limited amount of data in each batch, the true sparse structure of $\delta \textbf{f}$ cannot be accurately identified. It is often the case that $\mathbf{C}$ contains multiple entries that are nonzero due to the sampling error from the short time series. This means if directly imposing the sparsity into $\boldsymbol{\Xi}_r$ based on the CEM $\mathbf{C}$ computed from such a short single-batch data may lead to an incorrect residual system $\boldsymbol{\Xi}_r\boldsymbol{\Phi}$. To resolve such a sampling problem, we introduce aggregated CEM $\mathbf{C}^+(K)$ based on the average of the causation entropy from a series of data batches:
\begin{equation}
    \label{eq:C+}
    \mathbf{C}^+(K) =
    \begin{cases}
        0 \quad \text{if $\frac{1}{K}\sum_{k=1}^{K}C^{(k)}_{\phi_n \rightarrow r_i|[\boldsymbol{\Phi} \setminus \phi_n]} \leq \overline{C}$}, \\
        1 \quad \text{if $\frac{1}{K}\sum_{k=1}^{K}C^{(k)}_{\phi_n \rightarrow r_i|[\boldsymbol{\Phi} \setminus \phi_n]} > \overline{C}$}.
    \end{cases}
\end{equation}
Denote by $D$ the number of batches with which the aggregated CEM has not changed, namely
\begin{equation}
\label{eq:C_criterion}
    \mathbf{C}^+(K) = \mathbf{C}^+(K-d), \quad \text{for all $d=1,2,...,D-1$}.
\end{equation}
In the CEBoosting algorithm, $D$ is a hyper-parameter and needs to be pre-determined. We further define a stable aggregated causation entropy matrix $\overline{\mathbf{C}^+}=\mathbf{C}^+(K^*)$ with the smallest $K^*$ that satisfies the criterion in \eqref{eq:C_criterion}.

 With this stable aggregated causation entropy matrix $\overline{\mathbf{C}^+}$, we then impose sparsity into $\boldsymbol{\Xi}_r$ by setting $\delta\xi_{in}=0$ when $\overline{\mathbf{C}^+}_{in}=0$ and extract a set of remaining coefficients $\boldsymbol{\Xi}_r =\{\delta\xi_{in}|\overline{\mathbf{C}^+}_{in}=1\}$.

Once the sparsity is imposed into $\boldsymbol{\Xi}_r$ for extracting a set of remaining coefficients, a model $\boldsymbol{\Xi}_r\boldsymbol{\Phi}$ in \eqref{eq:residual_system} can be calibrated based on the accumulated batches of data $\textbf{x}(t)$ for $t \in [t_{B_1},t_{B_{K^*+1}})$. Assume a discrete approximation of the continuous data with a fixed time step $\Delta t$ such that $t_{B_{k+1}}-t_{B_k}=M\Delta t$ for any $k$. The model calibration is performed by solving the following least squares problem:
\begin{equation}
\label{eq:optim_single}
    \argmin_{\boldsymbol{\Xi}_r} \sum_{m=1}^{MK^*}\|\textbf{r}(t_{B_1}+m\Delta t) - \boldsymbol{\Xi}_r\boldsymbol{\Phi}(t_{B_1}+(m-1)\Delta t) \|^2,
\end{equation}
where $\|\cdot\|$ denotes the vector norm in $\bR^p$. Note that the method also works with adaptive time steps, and the assumption of a fixed time step in \eqref{eq:optim_single} is for the simplicity of the illustration.

\begin{figure}[H]
\centerline{
\includegraphics[width=1\textwidth]{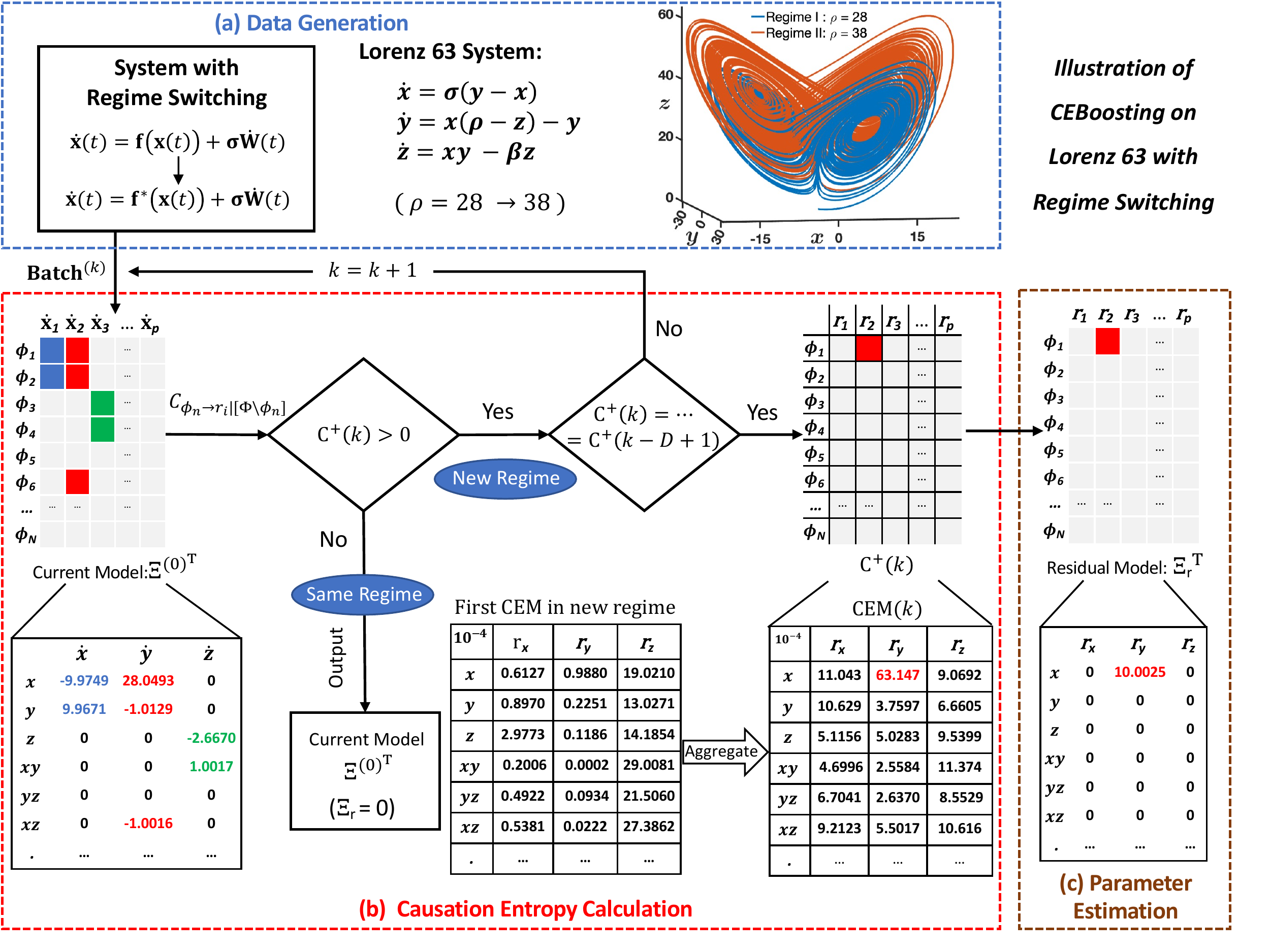}}
\caption{Schematic of CEBoosting algorithm (based on the Lorenz 63 model). Panel (a): data is generated from the Lorenz 63 system with regime switching. The parameter in Regime 1 is $\sigma=10, \beta=8/3, \rho=28$, while $\rho$ is changed to 38 in Regime 2. The index of incoming batch data $k$ starts with $1$. Panel (b): $\boldsymbol{\Xi}^{(0)}$ is the current model parameter matrix defined in \eqref{eq:discretized_system} with $\boldsymbol{\Phi}$ is the basis functions and $\dot{x}_i$ is the dynamic of state $x_i$. With $\boldsymbol{\Xi}^{(0)}$, the residual dynamics $r_i$ and the causation entropy between $r_i$ and $\boldsymbol{\Phi}$ can be calculated for each batch. $\text{CEM}(k)$ is the aggregated causation entropy from batch $1$ to $k$. $\mathbf{C}^+(k)$ is the binary matrix of $\text{CEM}(k)$ defined in \eqref{eq:C+} indicating the sparse structure of the residual model. $D$ is defined in \eqref{eq:C_criterion} indicating number that $\mathbf{C}^+(k)$ becomes stable, i.e., the pattern is consistent with the previous $D$ aggregated causation entropy matrix. Panel (c): with $\mathbf{C}^+(k)$ and data batches from the batch with new regime detected to the one with a stable $\mathbf{C}^+$, the residual model is calibrated by least square estimation.}
\label{fig:SchematicDiagram}
\end{figure}

After the sparse parameter matrix $\boldsymbol{\Xi}_r$ for residual dynamics model \eqref{eq:residual_system} is obtained, the current model $\boldsymbol{\Xi}^{(0)}$ is updated by adding the information from the residual dynamics $\boldsymbol{\Xi}_r$. The CEBoosting algorithm repeats the above procedure when a new batch of data arrives. A schematic illustration of Lorenz 63 system with regime switching is displayed in Fig.~\ref{fig:SchematicDiagram}, and more details of the CEBoosting algorithm can be found in Algorithm~\ref{alg:CEBoosting} presented in \ref{sec:algorithm}. It should be noted that the notations adopted in this section assume the regime switching time $t_s \in [t_{B_1},t_{B_2})$ for the simplicity of the illustration. In practice, regime switching can happen at $t_s \in [t_{B_k},t_{B_{k+1}})$ with $k>1$, for which Algorithm~\ref{alg:CEBoosting} presents the detailed procedures of detecting regime switching, aggregating causation entropy matrix, identifying a sparse model structure and then fitting the model parameters.

\subsection{Comparison with SINDy}

In recent years, SINDy has been a popular framework for sparsely identifying nonlinear dynamics from data \citep{brunton2016discovering}. SINDy exploits an iterative thresholding regularization to determine the model structure and estimate the model parameters simultaneously. The iterative thresholding regularization guarantees the parsimonious model structure. Using abundant training data (e.g., long time-series), the original SINDy method was designed for offline model identification. To better work with a limited amount of data, a recent extension of SINDy~\citep{fasel2022ensemble} leveraged the statistical approach of bagging and achieved a more robust learning performance.

Compared to the SINDy methods, the CEBoosting algorithm is mainly designed with a low computational cost for robust online learning with limited sequential data. The key difference from the SINDy methods is that imposing sparsity is decoupled from the parameter estimation in the CEBoosting algorithm. Specifically, causation entropy is utilized only to determine a parsimonious model structure without dealing with the model calibration. The concept of bagging is further introduced with sequential batches of data to ensure robust estimation of causation entropy before imposing a parsimonious model structure. It results in utilizing the minimum data to determine the model structure. With a robust estimation of the parsimonious model structure, the CEBoosting algorithm exploits a simple least square estimate to solve the parameter values, a quadratic optimization problem with a closed analytic solution. This avoids repeatedly estimating the parameters in the LASSO-type or thresholding-based regression approaches when determining the model structure by examining each batch data. In online learning, the CEBoosting algorithm, by design, can also avoid accidentally removing essential basis functions due to the incorrectly calibrated model based on a limited amount of data.

\section{Numerical Experiments for Systems with Moderately Large Dimensions, Partial Observations, and Extreme Events}
\label{sec:results}

In this section, the performance of the CEBoosting method is demonstrated utilizing four different chaotic or turbulent systems, which are models that mimic crucial features in many science or engineering disciplines. The presentation starts with the three-dimensional Lorenz 63 system, which is a classical chaotic system. This test is utilized as a proof of concept to demonstrate the detailed steps of the method in the online identification of the non-linear dynamics with regime switching. To illustrate the efficiency of the algorithm in capturing the regime switching behavior in a relatively high-dimensional case, the forty-dimensional Lorenz 96 system is utilized as a second test, where the localization technique is incorporated to mitigate the curse of dimensionality. As strongly non-Gaussian statistics, intermittency and extreme events appear in many climate, atmosphere and ocean science problems, a multi-mode layered topographic model that captures these crucial turbulent features is adopted as the next test model. The last test case aims to deal with a more realistic scenario where only a subset of the state variables is observed. Data assimilation is, therefore, essential in such a partial observational case. A stochastic parameterized extended Kalman filter (SPEKF) model is used to estimate and simulate the hidden system state variables. Below, assuming the starting model and its parameters are available for all the numerical examples. The goal is to learn the regime switching and the corresponding residual model that adjusts the original model to a new regime based on online sequential data.

For all four systems, the numerical simulation time step size is chosen as $\diff t=0.001$. The hyper-parameter $D$ in \eqref{eq:C_criterion} is chosen as 4 for the numerical examples of this work.

\subsection{The Lorenz 63 System: A Classical Chaotic System}
The Lorenz 63 (L63) system is proposed by Lorenz in 1963 \citep{lorenz1963deterministic}. It is a simplified mathematical model for atmospheric convection. The equations relate the properties of a two-dimensional fluid layer uniformly warmed from below and cooled from above. In particular, the equations describe the rate of change of three quantities concerning time: $x$ is proportional to the rate of convection, $y$  to the horizontal temperature variation, and $z$  to the vertical temperature variation. The constants $\sigma$, $\rho$, and $\beta$ are system parameters proportional to the Prandtl number, Rayleigh number, and certain physical dimensions of the layer itself \citep{sparrow2012lorenz}. The L63 model is also widely used as a simplified model for lasers, dynamos, thermosyphons, brushless DC motors, electric circuits, chemical reactions, and forward osmosis \citep{haken1975analogy, knobloch1981chaos, gorman1986nonlinear, hemati1994strange, cuomo1993circuit, poland1993cooperative, tzenov2014strange}.
The governing equation of the Lorenz 63 system is as follows,
\begin{equation}
\begin{aligned}
&\frac{\diff x}{\diff t} = \sigma (y-x), \\
&\frac{\diff y}{\diff t} = x(\rho -z)  - y,\\
&\frac{\diff z}{\diff t} = xy - \beta z.
\end{aligned}
\label{eq:lorenz63}
\end{equation}
The standard parameters $\sigma=10$, $\beta=8/3$, and $\rho = 28$ that create the butterfly profile are utilized as the starting regime. The new regime takes a different value of the parameter $\rho=38$. The sudden change of the parameter from $\rho = 28$ to $\rho = 38$ occurs at $t=100$. The goal is to (i) detect the regime switching and (ii) learn the residual model that adjusts the original system to the new one.

\begin{figure}[H]
\includegraphics[scale=0.4]{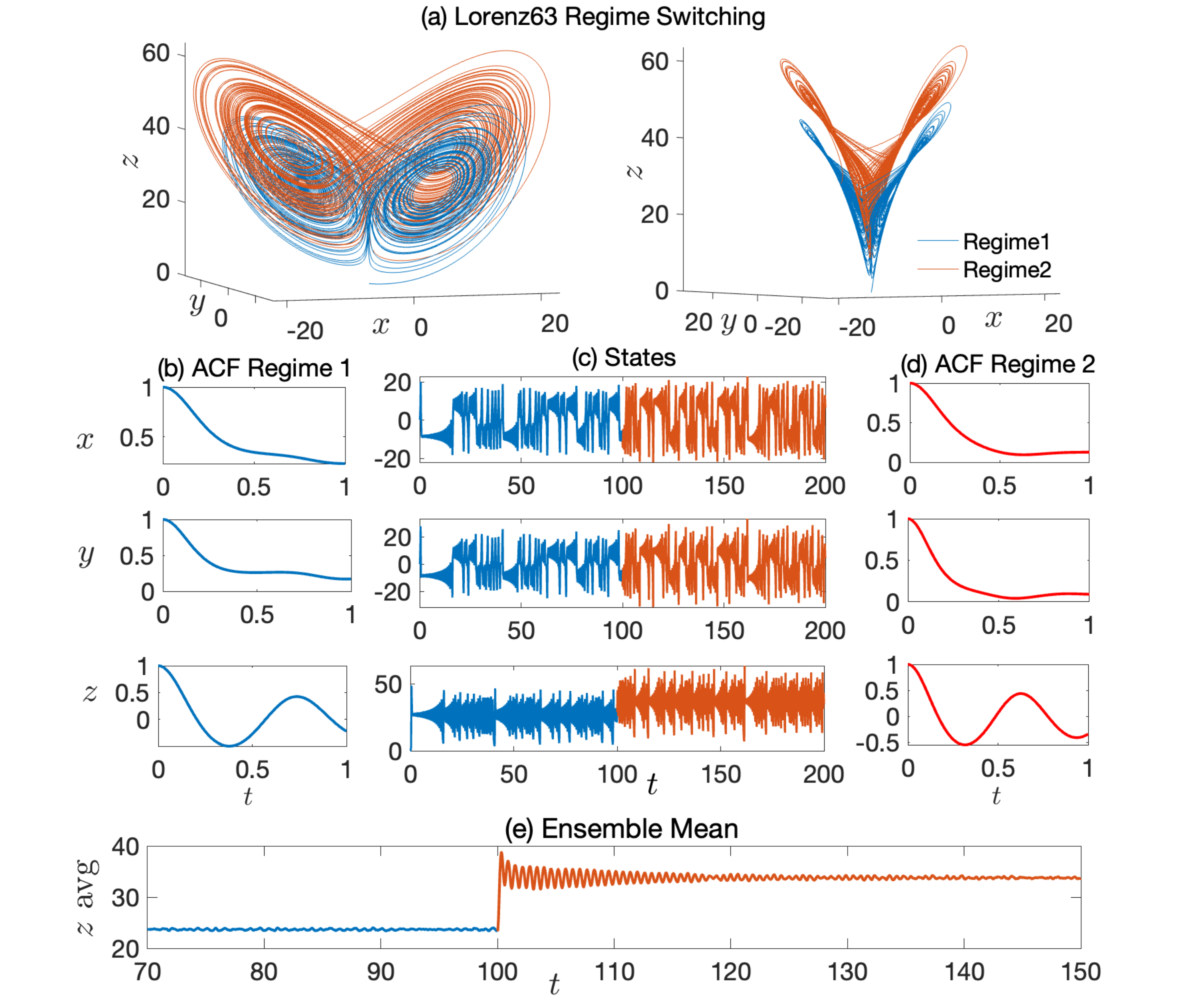}
\centering
\caption{Lorenz 63 system with regime switching. (a): trajectories of the original system and the new one in phase space. (c): time series of system state variables. (b) and (d): autocorrelation function (ACF) of the original system and the new regime. (e): ensemble mean of variable $z$ before and after regime switching at $t=100$.}
\label{fig:lorenz63}
\end{figure}

Figure~\ref{fig:lorenz63}(a) compares the trajectories in the phase space between the original (Regime 1) and new (Regime 2) regimes. It can be seen that the two systems are on different manifolds in the phase space. Figure~\ref{fig:lorenz63}(c) shows the time series of the three state variables. With the regime switching at $t=100$, it can be seen in Fig.~\ref{fig:lorenz63}(c) that the patterns of time series change accordingly, especially for the variable $z$ that demonstrates a shift of its mean value. The autocorrelation function (ACF) of each state variable for the original system is presented in Fig.~\ref{fig:lorenz63}(b), which shows a rapid decay of correlation within one time unit, except for the variable $z$ that has some oscillations in its ACF. Similar behavior of the ACF is observed in Fig.~\ref{fig:lorenz63}(d) after the regime switching. Finally, Figure~\ref{fig:lorenz63}(e) shows the ensemble mean of the state variable $z$ as a function of time from an independent simulation with 5000 ensemble members, which reveals that the transition time of the regime switching is about $20$ time units. This also indicates that the transition time depends on the property of the transient feature and is very different from the decorrelation time. Nevertheless, the decorrelation time of the original system provides a natural way to determine the batch size. Thus, the batch size is chosen as one time unit here.

The CEBoosting algorithm is employed to detect the regime switching and identify the sparse structure of the residual model. The candidate basis functions include all the linear and quadratic nonlinear functions: $\{x, y, z, xy, xz, yz, x^2, y^2, z^2\}$. With the identified sparse structure, the coefficients of the residual model are then determined via the least square estimation.
To detect the regime switching and identify the sparse structure of the residual model, the aggregated causation entropy matrix is gradually updated until its structure gets stable, i.e., the causation entropy values of some candidate basis functions keep being significantly greater than others. The stable causation entropy matrix structure is obtained after 6 time units and is presented in Table~\ref{tab:lorenz63_CEM}, where only one entry has a significant value. It is worthwhile to highlight that the total time units (6 units) to discover the regime switching and determine the model parameters in the new regime is shorter than the transition time (20 units). In other words, utilizing the information from the transient period is sufficient for the CEBoosting algorithm to determine the model response to the regime switching in this test case.

\begin{table}[H]
\begin{center}
\caption{The L63 model - The causation entropy matrix (CEM) after 6 time units. The pattern of CEM becomes stable and does not change with incorporating more batch data. The entry with a significant value of the causation entropy is highlighted using the bold font. According to the pattern of this CEM, a residual model will be built.}
\label{tab:lorenz63_CEM}
\begin{tabular}{ |c|c|c|c|c|c|c|c|c|c| }
 \hline
$10^{-4}$& $x$ & $y$ & $z$ & $xy$ & $xz$ & $yz$ & $x^2$ & $y^2$ & $z^2$ \\
 \hline
 $\dot{x}$ &11.0437  & 10.6292  &  5.1156 &   4.6996  &  6.7041 &   9.2123&    8.9308&    3.0866  &  4.0888 \\
 \hline
 $\dot{y}$ & \textbf{63.1474}&   3.7597 &   5.0283  &  2.5584   & 2.6370    &5.5017   & 4.3689   & 2.4472 &   4.4099  \\
 \hline
 $\dot{z}$ & 9.0692   & 6.6605 &   9.5399 &  11.3740  &  8.5529 &  10.6165&   10.8299  & 10.8725  & 10.5916  \\
 \hline
\end{tabular}
\end{center}
\end{table}

Finally, based on the sparse structure identified in Table~\ref{tab:lorenz63_CEM}, a residual model can be calibrated via the least square estimation based on the residuals of the original model. The residual model is then added to the original model as a correction term. The coefficients of the corrected model for the new regime are listed in Table~\ref{tab:lorenz63_model}. It can be seen that the corrected model successfully captures the crucial parameter value $\rho=38$ that induces the new regime.

\begin{table}[H]
\begin{center}
\caption{The L63 model - Updated model parameters for the new regime.}
\label{tab:lorenz63_model}
\begin{tabular}{ |c|c|c|c|c|c|c|c|c|c| }
 \hline
  & $x$ & $y$ & $z$ & $xy$ & $xz$ & $yz$ & $x^2$ & $y^2$ & $z^2$ \\
 \hline
 $\dot{x}$ &   -9.9749  & 9.9671 & 0 & 0 & 0 & 0 & 0 & 0 & 0 \\
 \hline
 $\dot{y}$ & \textbf{38.0357}    & -1.0129   & 0 & 0 & 0 & -1.0016 & 0 & 0 & 0  \\
 \hline
 $\dot{z}$ & 0 & 0 &  -2.6670  & 1.0017  & 0 & 0 & 0 & 0 & 0  \\
 \hline
\end{tabular}
\end{center}
\end{table}

\subsection{The Lorenz 96 System: Strategy for Applying CEBoosting to Relatively High-Dimensional Systems}

\sloppy The Lorenz 96 (L96) system is a classical chaotic to turbulent model \citep{lorenz1996predictability}. It can be regarded as a coarse discretization of atmospheric flow on a latitude circle with complicated wave-like and chaotic behavior. It is also widely used as a testbed for data assimilation, prediction, and uncertainty quantification \citep{wilks2005effects, lee2017multiscale, arnold2013stochastic, chen2017beating}. The L96 model reads:
\begin{equation}
\frac{\diff x}{\diff t} = (x_{j+1} - x_{j-2})x_{j-1} - x_j + F,\qquad
j = 1, 2, ..., J,
\label{eq:lorenz96}
\end{equation}
with periodic boundary conditions $x_{-1} = x_{J-1}, x_0 = x_J, x_{J+1} = x_1$ and $J=40$ is adopted.
The L96 system is utilized here to demonstrate that the CEBoosting algorithm, combined with localization techniques, can detect regime switching for relatively high-dimensional systems. The motivation here is that the number of candidate functions quickly increases with the dimension of the system if all possible combinations of state variables up to a particular order are included. As a result, the computational cost becomes unaffordable without suitable treatment for such a curse of dimensionality. To this end, the idea of localization is exploited here. Localization means the dynamics of each state variable only depend on the nearby ones. In fact, the advection, diffusion, and dispersion are all local operators \citep{majda2003introduction, vallis2017atmospheric}. Similarly, the localization applies to many parameterization problems for the subgrid scales, which also depend only on the nearby corresponding large-scale state variables \citep{grabowski2004improved, gagne2020machine, chattopadhyay2020data}. The idea of localization is also widely utilized in data assimilation and prediction \citep{bergemann2010localization, anderson2007exploring, janjic2011domain}. Localization reduces the number of candidate functions for the dynamics of each $x_j$ by including only the terms that represent local interactions with $\diff x_j/\diff t$. Note that the total number of candidate functions in the entire library can remain large, but only a relatively small number of functions will be examined for the causal relationship to the dynamics of each $x_j$.

\begin{figure}[H]
\includegraphics[scale=0.4]{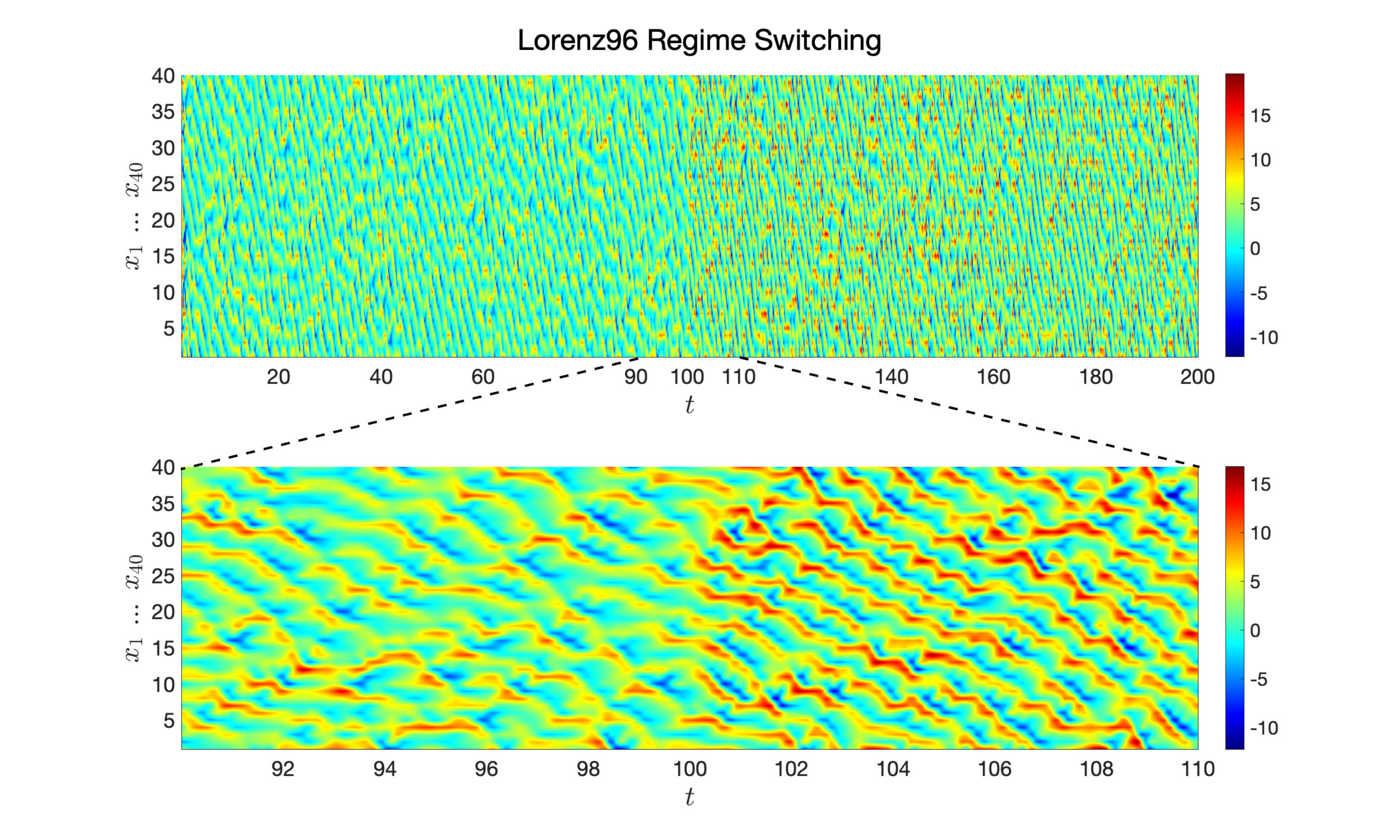}
\centering
\caption{The Lorenz 96 system with regime switching. Top: forty-dimensional system state evolving at the time interval $[0, 200]$ with a regime switching at $t=100$. Regime 1 is chaotic with $F=8$, and regime 2 is turbulent with $F=16$ and the linear terms $-1.5x_j$. Bottom: zoom-in view of the regime switching at the time interval $[90, 110]$.}
\label{fig:l96_pcolor}
\end{figure}

\begin{figure}[H]
\begin{center}
\includegraphics[scale=0.4]{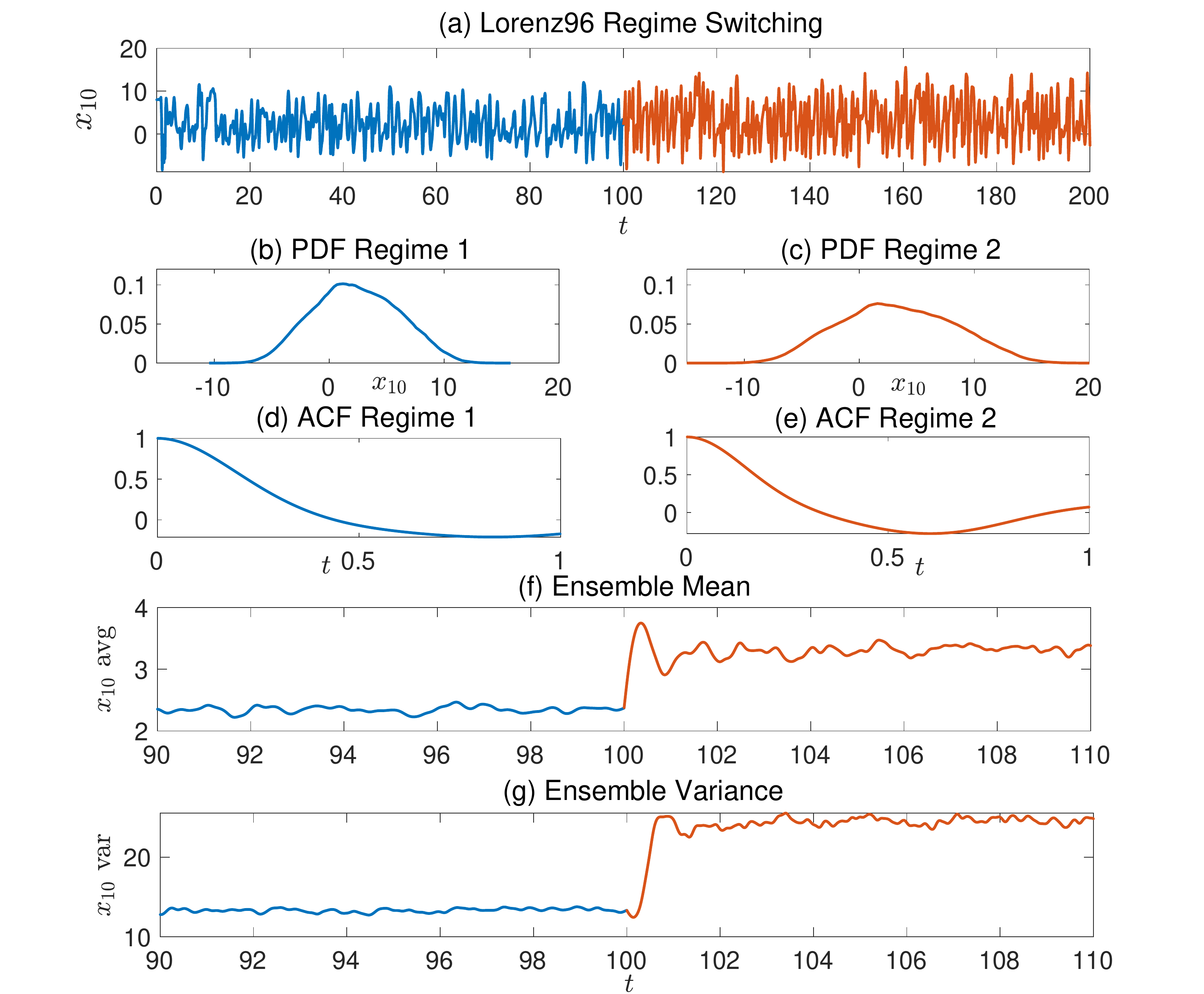}
\end{center}
\caption{Statistical properties of the state variable $x_{10}$ of the Lorenz 96 system with regime switching. Panel (a): time series of $x_{10}$ in regime 1 and regime 2. Panels (b) and (c): the PDFs of $x_{10}$ in both regimes. Panels (d) and (e): the ACFs of $x_{10}$ in both regimes. Panels (f) and (g): the ensemble mean and the ensemble variance of $x_{10}$ before and after the regime switching at $t=100$.}
\label{fig:l96_synthesis}
\end{figure}

The original system (Regime 1) has a forcing coefficient $F = 8$, which makes the system have a strongly chaotic behavior. After the regime switching at $t=100$, the new system (Regime 2) takes the forcing value $F=16$. Meanwhile, the coefficient of $x_j$, representing the damping effect, changes from $-1$ to $-1.5$. This leads to a fully turbulent regime \citep{majda2012filtering}. Some weak coherent structures can still be observed in the strongly chaotic regime, but they disappear in the fully turbulent one. See Figure~\ref{fig:l96_pcolor} for the two regimes. Figure~\ref{fig:l96_synthesis} shows the detailed statistical properties of the state variable $x_{10}$. Note that the model has homogeneous dynamics, so the statistical properties for different state variables are the same.
Figure~\ref{fig:l96_synthesis}(a) illustrates how the time series of $x_{10}$ change with the regime switching. From the strongly chaotic to the fully turbulent regime, the variance of all the state variables becomes more extensive, and the decorrelation time becomes shorter. The batch size of 1 time unit is chosen here, which is again of the same order as the decorrelation time. The total transition time is about 1 to 2 units, according to Figure~\ref{fig:l96_synthesis}(b).

A set of candidate basis functions $\{x_j, x_j^2, x_jx_{j-1}, x_jx_{j+1}, x_jx_{j-2}, x_jx_{j+2}\}$ is built by exploiting polynomials up to the second order for each state $x_j,\ j=1,2, ...J$. Note that the residual dynamics of each state $x_j$ is assumed to have contributions from the basis functions that consist of only its adjacent states.
The CEBoosting algorithm identifies the sparse structure of the residual model using about two batches of data (i.e., the causation entropy matrix pattern does not change after two time units). Table~\ref{tab:lorenz96_CEM} shows the CEM based on utilizing two batches of data, which is about the same length as the estimated transition time of Lorenz 96 system (see Figure~\ref{fig:l96_synthesis}(f-g)). It is seen that the causation entropy entries associated with the actual residual terms are at least one order more significant than others. Therefore, the CEBoosting algorithm successfully identifies the correct sparse structure of the residual model based on data within the time interval of the transient period.

\begin{table}[H]
\begin{center}
\caption{The L96 model  - The CEM after 2 time units. The entries with significant values of the causation entropy are highlighted using the bold font. The pattern of CEM becomes stable and does not change with incorporating more data batches. According to the pattern of this CEM, a residual model will be built.}
\label{tab:lorenz96_CEM}
\begin{adjustbox}{max width=1.\textwidth,center}
\begin{tabular}{ |c|c|c|c|c|c|c|c|c|c|c|c|c|c| }
 \hline
  $10^{-4}$& $x_1$ & $x_1^2$ & $x_2x_{40}$ & $x_{39}x_{40}$ & $x_2$ & $x_2^2$ & $x_3x_1$ & $x_{40}x_1$ & ...& $x_{40}$ & $x_{40}^2$ & $x_1x_{39}$ & $x_{38}x_{39}$ \\
 \hline
 $\dot{x_1}$ &\textbf{436.0597}& 4.2015&2.2964&1.4823&0.0976 &3.4716&1.3792  &0.9251&...&0.1869&3.0008&0.7236&8.0802 \\
 \hline
 $\dot{x_2}$    & 8.4920&    9.8888&    6.0428 &   8.0793&  \textbf{19.8015}&    8.3095  & 10.1492 &   9.6844 & ... &   0  &  0  &  0 & 0
 \\
 \hline
 ... & ... &...&...&...&...&...&...&...&...&...&...&...&...\\
 \hline
 $\dot{x_{40}}$   &  2.7772  &  3.1363  &  4.5714 &   8.4405    &9.6327  &  5.3498 &   4.8701 &  7.6239 & ...&  \textbf{200.3710}  &  2.5200  &  7.6225 &   2.7748
 \\
 \hline
\end{tabular}
\end{adjustbox}
\end{center}
\end{table}

After obtaining the sparse structure according to the CEM in Table~\ref{tab:lorenz96_CEM}, the residual dynamics are calibrated via a linear combination of candidate basis functions, and the coefficients for the linear combination is obtained by the least square estimation. The L96 model in the new regime is then updated by adding this residual model to the original model. The coefficients of the updated model are summarized in Table~\ref{tab:lorenz96_model}. The identified model shows a good agreement with the true system of the new regime with $F=16$ and the linear term $-1.5x_j$.

\begin{table}[H]
\begin{center}
\caption{The L96 model - Coefficients of the updated model for the new regime. The entries with bold font indicate the parameters that contribute to the regime switching.}
\label{tab:lorenz96_model}
\begin{adjustbox}{max width=1.\textwidth,center}
\begin{tabular}{ |c|c|c|c|c|c|c|c|c|c|c|c|c|c|c| }
 \hline
  & $1$ & $x_1$ & $x_1^2$ & $x_2x_{40}$ & $x_{39}x_{40}$ & $x_2$ & $x_2^2$ & $x_3x_1$ & $x_{40}x_1$ & ...& $x_{40}$ & $x_{40}^2$ & $x_1x_{39}$ & $x_{38}x_{39}$ \\
 \hline
 $\dot{x_1}$ &\textbf{16.0248}&\textbf{-1.4938}&0&1.0000&-1.0016 &0&0&0&0&...&0&0&0&0 \\
 \hline
 $\dot{x_2}$ &\textbf{16.0490}&0&0&0&0&\textbf{-1.5077}&0&1.0000&-1.0000&...&0&0&0&0  \\
 \hline
 ... & ... &...&...&...&...&...&...&...&...&...&...&...&...&...\\
 \hline
 $\dot{x_{40}}$ &\textbf{16.0027}&0&0&0&0&0&0&0&0&...&\textbf{-1.5017}&0&1.0000&-1.0000 \\
 \hline
\end{tabular}
 \end{adjustbox}

\end{center}
\end{table}

\subsection{The Topographic Model: Dynamical System with Intermittency and Extreme Events}
The topographic model is an ideal model to study the complex nonlinear interaction of the large-scale and the small-scale flow and the role of the topography \citep{majda2006nonlinear, majda2003systematic}. The topographic model can generate intermittency and extreme events. The corresponding PDF is often highly non-Gaussian with heavy tails. Therefore, it provides a challenging test case for detecting the regime switching behavior.
Here the small-scale flow is given in terms of the stream function $\psi$. The large-scale velocity field only has the zonal component $u(t)$, and the topography is given by the function $h$. The parameter $\beta > 0$ is the contribution from the beta-plane effect. A common simplified version of the topographic model assumes a layered topography along the $y$ direction, which means $\psi$ is only a function of $y$. In addition, the model contains only the leading two Fourier wavenumbers of the stream function (with $k = \pm1$ and $\pm2$). With these simplifications, the resulting model is reduced to a $5$
dimensional system containing $u$, $\psi_{\pm1}$ and $\psi_{\pm2}$. For the simplicity of notation, a change of variables defines the new state variables $v_1,\ldots,v_4$, which are linked with $\psi_{\pm1}$ and $\psi_{\pm2}$ via
\begin{equation*}
  \psi_1 = \frac{1}{2\sqrt{2}}\left((v_2-v_1)-(v_2+v_1)i\right)\qquad\mbox{and}\qquad \psi_2 = \frac{1}{2\sqrt{2}}\left((v_4-v_3)-(v_4+v_3)i\right)
\end{equation*}
where $i$ is the imaginary unit.
Similarly, $\omega_1$ and $\omega_3$ are the two new variables standing for the Fourier coefficients of the topographic effect from $h$. The model reads:
\begin{equation}
\begin{aligned}
&\frac{\diff v_1}{\diff t} = -d_{v_1}v_1 - \beta v_2 + v_2u - 2\omega_1 u + \sigma_{v_1} \dot{W_1} \\
&\frac{\diff v_2}{\diff t} = -d_{v_2}v_2 - \beta v_1 - v_1u + \sigma_{v_2} \dot{W_2}\\
&\frac{\diff v_3}{\diff t} = -d_{v_3}v_3 - \omega_3u - \frac{\beta}{2}v_4 + 2v_4u + \sigma_{v_3} \dot{W_3} \\
&\frac{\diff v_4}{\diff t} = -d_{v_4}v_4 - \frac{\beta}{2}v_3 - 2v_3u + \sigma_{v_4} \dot{W_4} \\
&\frac{\diff u}{\diff t} = -d_uu + \omega_1v_1 + 2\omega_3v_3 + \sigma_u \dot{W_u}
\end{aligned}
\label{eq:topographical}
\end{equation}
The original system (Regime 1) takes the parameters $d_{v_1} = d_{v_2} = d_{v_3} = d_{v_4} = d_{v_5} = 0.005$, $\beta = 1$,  $\sigma_{v_1} = \sigma_{v_2}= \sigma_{v_3}= \sigma_{v_4} = 1/(20\sqrt 2)$, $\sigma_u = 1/\sqrt 2$, $\omega_1 = \sqrt{2}/2$ and $\omega_3 = \sqrt{2}/4$. The new model (Regime 2) utilizes different parameters for the topographic effect with $\omega_1 = \omega_3 = 3\sqrt{2}/2$. It is worthwhile to note that the noise level in the dynamics of the zonal flow $u$ is much larger than those in the dynamics of $v_i$ in \eqref{eq:topographical}. This is a typical situation as the $u$ is the only mode that explicitly describes the zonal feature of the flow. A large noise is taken to mimic the unresolved dynamical features.

\begin{figure}[H]
\includegraphics[scale=0.4]{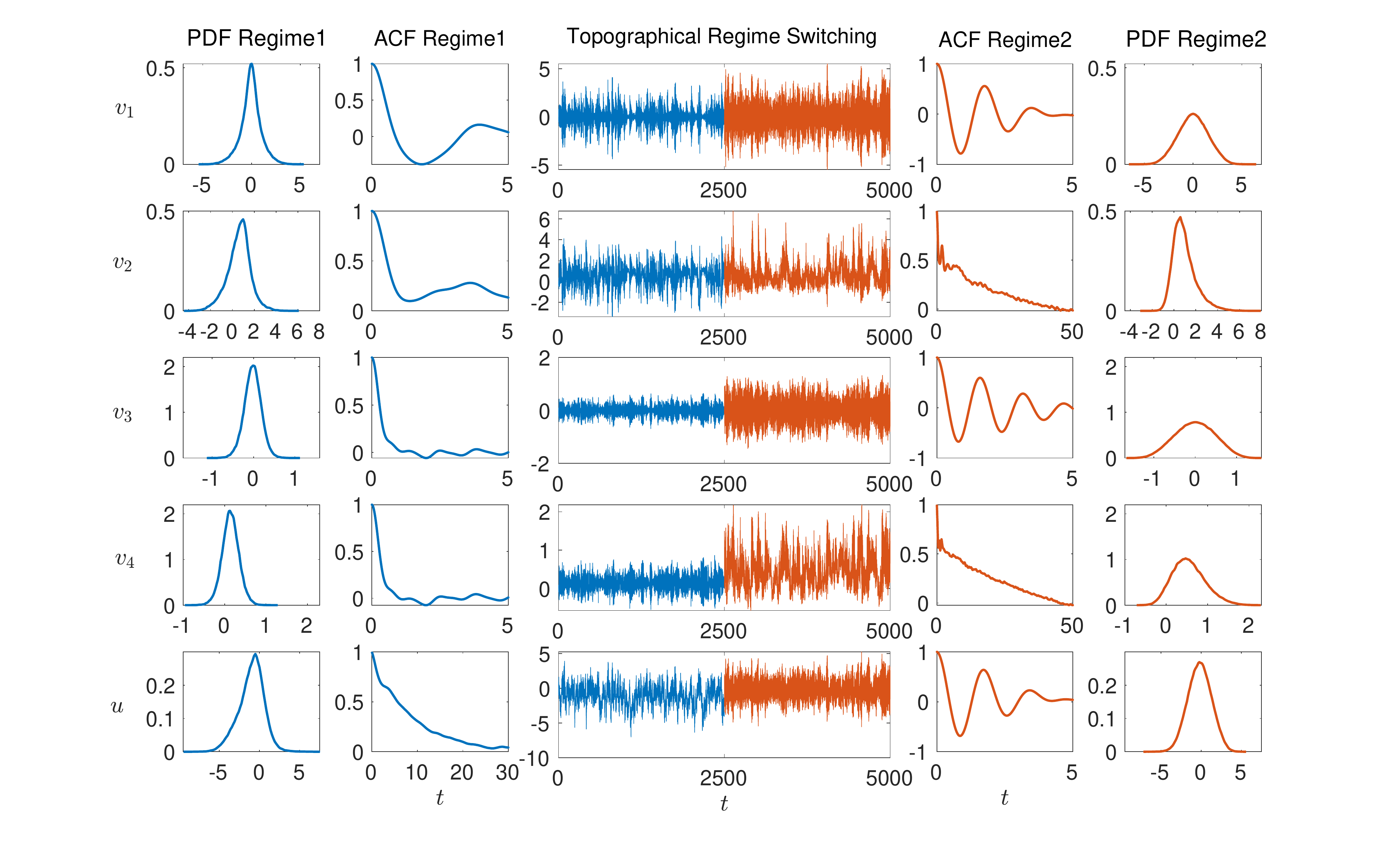}
\centering
\caption{Topographic model with regime switching. Middle column: time series of each state in the time interval [0, 5000] with regime switching at $t=2500$. 1st and 5th columns: probability density function (PDF) of each state in both regimes. 2nd and 4th columns: autocorrelation function (ACF) of each state in both regimes.}
\label{fig:topo_synthesis}
\end{figure}

Figure \ref{fig:topo_synthesis} displays the regime switching and statistical properties of all five variables in the topographic model. The pattern and amplitude of each time series demonstrate a noticeable change after the regime switching happens at $t=2500$. Note that the long time series of both regimes are utilized here, ensuring that the associated statistics are computed with a sufficiently large number of sample points. The PDFs behave like Laplace distributions with heavy tails, indicating many extreme events in the model simulation. The ACF of $u$ in the original model (regime 1) decays very slowly, taking about 30 time units until the ACF approaches zero. Because of this, a batch size of 30 time units is utilized in the CEBoosting algorithm. Due to the strengthening of the topographic effect after the regime switching, the time series of $v_2$ and $v_4$ in the new model (regime 2) display multiscale features, where the ACFs have a quick decay at the beginning but then relax to zero slowly after $50$ time units.

Similar to the previous test models, the library contains the polynomials up to the second-order. In other words, the following twenty basis functions are employed to build the library: $\{v_1, v_2, v_3, v_4, u,v_1^2, v_2^2, v_3^2, v_4^2, u^2, v_1v_2, v_1v_3, v_1v_4, v_1u, v_2v_3, v_2v_4, v_2u, v_3v_4, v_3u, v_4u\}$.

As the noise levels of state variables are significantly different from each other, the magnitudes of causation entropy for different residual dynamics are not the same either. Therefore, the selection of the candidate functions is based on the causation entropy values for the dynamics of each state variable separately. To determine the entire residual model, 60 time units is used. See Table~\ref{tab:topographical_CEM_60}.

\begin{table}[H]
\begin{center}
\caption{The topographic model - The CEM after 60 time units. The pattern of CEM does not change with incorporating more batch data. The entries with significant values of the causation entropy are highlighted using the bold font.}
\label{tab:topographical_CEM_60}
\begin{tabular}{ |c|c|c|c|c|c|c|c|c|c|c| }
 \hline
  $10^{-3}$& $v_1$ & $v_2$ & $v_3$ & $v_4$ & $u$ & $v_1u$ & $v_2u$ & $v_3u$ & $v_4u$ & ... \\
 \hline
 $\dot{v_1}$ &0.0115 & 0.0227 & 0.0165 & 0.0614 &\textbf{120.9879} &0.0216 & 0.0365&0.0407&0.0249&... \\
 \hline
 $\dot{v_2}$& 0.0124 &0.0327&0.0066&0.0150 & 0.0313& 0.0002  &  0.0090 & 0.0039 & 0.0284&...  \\
 \hline
 $\dot{v_3}$ &0.0023&    0.0033&    0.0142&    0.0099&   \textbf{51.3949}   & 0.0537   &   0.0291  &  0.0572  & 0.0082&...\\
 \hline
 $\dot{v_4}$ &0.0359   & 0.0619 &   0.0150 &   0.0134&    0.0253   &0.0180  &    0.0127 &0.0099  & 0.0210&...\\
 \hline
 $\dot{u}$ & \textbf{0.1080}  &  0.0022 &   \textbf{0.0471}  &  0.0146 &   0.0065    & 0.0068 &  0.0097 & 0.0049 &  0.0063 & ... \\
 \hline
 \end{tabular}
\end{center}
\end{table}

Based on the sparse model structure identified in Table~\ref{tab:topographical_CEM_60}, the residual model is calibrated via the least square estimation. The coefficients of the calibrated model for the new regime are summarized in Table~\ref{tab:topographical_model}. It can be seen that the corrected model successfully updates the changed parameters $\omega_1$ and $\omega_3$ in the new regime. The resulting model can reproduce the strong non-Gaussian features with intermittency and extreme events.

\begin{table}[H]
\begin{center}
\caption{The topographic model - Updated model for the new regime}
\label{tab:topographical_model}
\begin{tabular}{ |c|c|c|c|c|c|c|c|c|c|c| }
 \hline
  & $v_1$ & $v_2$ & $v_3$ & $v_4$ & $u$ & $v_1u$ & $v_2u$ & $v_3u$ & $v_4u$ & ... \\
 \hline
 $\dot{v_1}$ &-0.0496 & -1.0008&0&0&\textbf{-4.2463} &0&0.9998&0&0&... \\
 \hline
 $\dot{v_2}$ &0.9988& -0.0501&0&0&0&-1.0009&0&0&0&...  \\
 \hline
 $\dot{v_3}$ &0&0&-0.0554& -0.4968&\textbf{-2.1200}  &0&0&0&2.0015&...\\
 \hline
 $\dot{v_4}$ &0&0&0.4998&-0.0523&0&0&0&-1.9996 &0&...\\
 \hline
 $\dot{u}$ &\textbf{2.1123}&0& \textbf{4.5713}&0&-0.0599&0&0&0&0&... \\
 \hline
 \end{tabular}
\end{center}
\end{table}

\begin{figure}[H]
\includegraphics[scale=0.4]{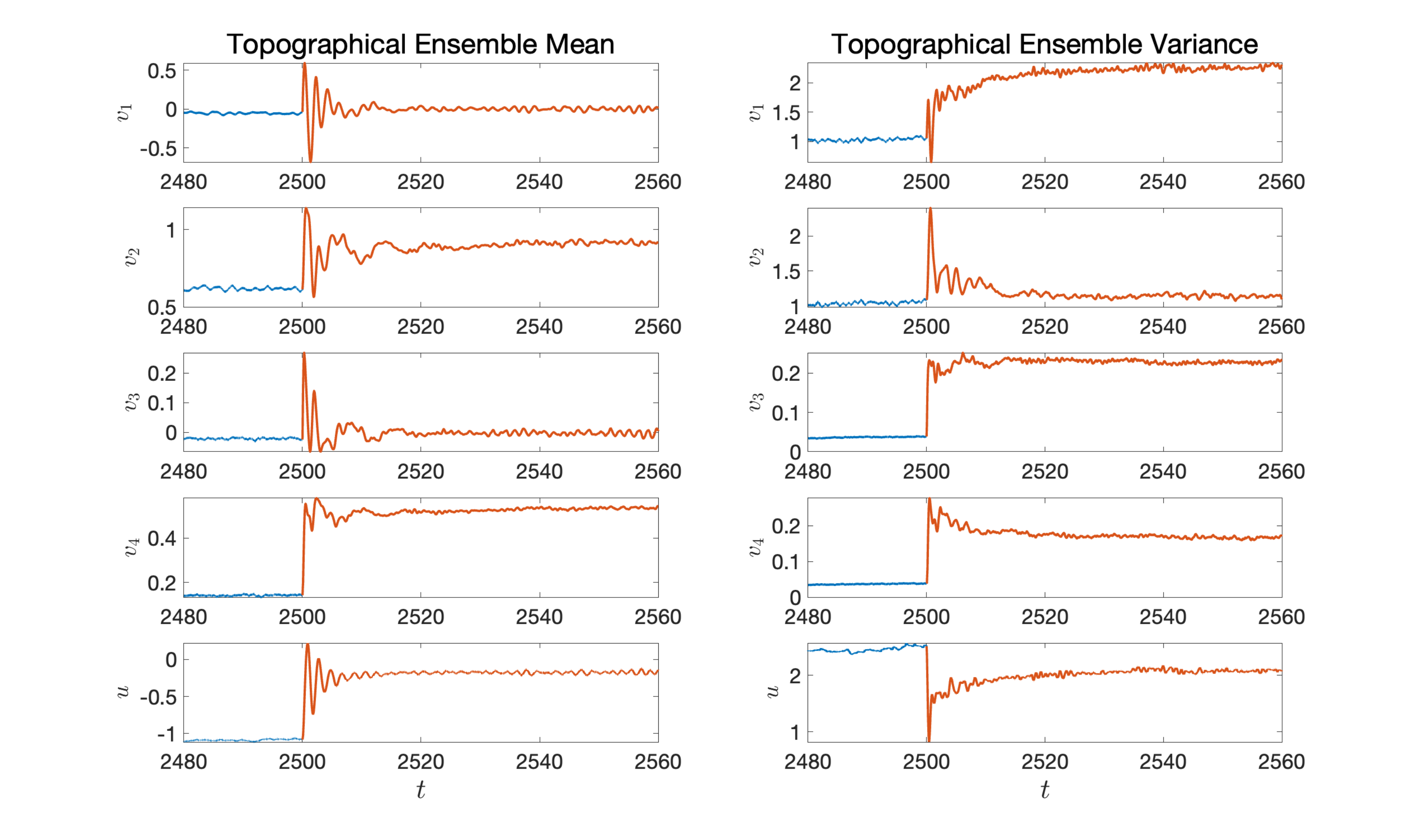}
\centering
\caption{Ensemble mean and variance of each state variable in the topographic model. Regime switching happens at $t = 2500$.}
\label{fig:topo_ensemble}
\end{figure}

It is worthwhile to remark that, as the noise levels in the dynamics of $v_i$ are lower than that in $u$, a shorter time series (and shorter batch length) with in total of only 15 time units can be utilized to reach a stable CEM for the $v_i$ components. See Table~\ref{tab:topographical_CEM_15}. Figure~\ref{fig:topo_ensemble} shows the time evolution of the ensemble mean and the ensemble variance of the topographic model, including the time instant of regime switching. The results here can be utilized to infer the transition time, which is about 20 time units. Therefore, in this topographic model with a large noise in the $u$ dynamics, the stable pattern of all $v_i$ variables can be identified by the CEBoosting algorithm within the transition period. Yet, 60 time units of data are needed for identifying the model structure associated with the variable $u$, mainly due to the larger noises of $u$.

\begin{table}[H]
\begin{center}
\caption{The topographic model - The CEM after 15 time units. The pattern of all $v$ variables is stable. The entries with significant values of the causation entropy are highlighted using the bold font.}
\label{tab:topographical_CEM_15}
\begin{tabular}{ |c|c|c|c|c|c|c|c|c|c|c| }
 \hline
$10^{-3}$ & $v_1$ & $v_2$ & $v_3$ & $v_4$ & $u$ & $v_1u$ & $v_2u$ & $v_3u$ & $v_4u$ & ... \\
 \hline
 $\dot{v_1}$ &0.0182 &0.0233 &0.0322 &0.0679 &\textbf{22.665} &0.0244 &0022&0.1726&0.0285&... \\
 \hline
 $\dot{v_2}$ &0.0768 &0.0071 & 0.0002&0.0004  &  0.0879&  0.0582&0.0007&0.0171 &0.1059&...  \\
 \hline
 $\dot{v_3}$ &0.0187 &0.0006 & 0.1123& 0.0411 & \textbf{10.3880} & 0.0133&0.0024&0.0843 &0.0166&...\\
 \hline
 $\dot{v_4}$ &0.0087 &0.1908 &0.0171& 0.0023&0.000&0.0001&0.0520& 0.0015 &0.0016&...\\
 \hline
 $\dot{u}$ &0.0008 &0.0007 &0.0078 &0.0194   & 0.0167&0.0250 &0.0054 & 0.0223 &0.0107&... \\
 \hline
 \end{tabular}
\end{center}
\end{table}

\subsection{A Stochastic Parameterized Extended Kalman Filter (SPEKF) Model: Incorporating Data Assimilation into the CEBoosting Algorithm}

In many practical situations, the observations are only available for a subset of state variables, known as partial observations. It will be shown in the following that the CEBoosting algorithm can be naturally applied to the case with partial observations when data assimilation is appropriately incorporated. To illustrate the CEBoosting algorithm in the partial observational scenario, a simple yet practically useful nonlinear model is utilized as a testbed. The model is the so-called stochastic parameterized extended Kalman filter (SPEKF) model \citep{gershgorin2010improving, gershgorin2010test},
\begin{equation}
\begin{aligned}
&\frac{\diff u}{\diff t} = [-(\gamma + \hat{\gamma}) + i(\omega + \hat{\omega})]u  + (b+\hat{b}) + \sigma_u\dot{W}_u, \\
&\frac{\diff \gamma}{\diff t} = -d_{\gamma}\gamma + \sigma_\gamma \dot{W}_\gamma, \\
&\frac{\diff \omega}{\diff t} = -d_{\omega}\omega + \sigma_\omega \dot{W}_\omega, \\
&\frac{\diff b}{\diff t} = -d_{b}b + \sigma_b \dot{W}_b.
\end{aligned}
\label{eq:SPEKF}
\end{equation}
In this SPEKF model, $u(t)$ is a complex-valued state variable and is the only variable in the system to be observed. The observed variable $u(t)$ is driven by three hidden variables $\gamma(t), \omega(t)$ and $b(t)$. The parameters $d_\gamma, d_\omega, d_b$ are all positive, serving as damping factors. The parameters $\sigma_u$, $\sigma_\gamma$, $\sigma_\omega$ and $\sigma_b$ are noise coefficients, which are also positive. The white noises $\dot{W}_u$ and $\dot{W}_b$ are complex-valued while $\dot{W}_\gamma$ and $\dot{W}_\omega$ are real. The governing equations of $\gamma(t)$, $\omega(t)$ and $b(t)$ are Ornstein–Uhlenbeck (OU) processes \citep{gardiner2004handbook} with $\gamma$ and $\omega$ taking real values and $b$ taking a complex value. The three constants $\hat\gamma$, $\hat\omega$ and $\hat{b}$ in the dynamics of $u$ represent the mean damping, mean phase, and the mean forcing, respectively.

The SPEKF model \eqref{eq:SPEKF} has been widely used as an approximate model to describe a spectral mode of a complex turbulent system, especially in the context of data assimilation and ensemble prediction \citep{chen2023stochastic, majda2012filtering, gershgorin2010test}. Physically, the variable $u(t)$ represents one of the resolved modes (i.e., observable) in the turbulent signal, while the three hidden variables $\gamma(t), \omega(t)$ and $b(t)$  are surrogates for the nonlinear interaction between $u(t)$ and other unobserved modes in the original governing equation after applying the spectral decomposition. The idea of the SKEPF model is that the small or unresolved scale variables are stochastically parameterized by inexpensive linear and Gaussian processes, representing stochastic damping $\gamma(t)$, stochastic phase $\omega$ and stochastic forcing $b$. Despite the model error in using such Gaussian approximations for the original unresolved nonlinear dynamics, these Gaussian processes succeed in providing accurate statistical feedback from the unresolved scales to the resolved ones. Thus the intermittency and non-Gaussian features observed in the resolved variables can be accurately recovered. The statistics in the SPEKF model can also be solved with exact and analytic formulae, which allow an accurate and efficient estimation of the model states. The SPEKF type of model has been used for filtering multiscale turbulent dynamical systems \citep{majda2012filtering}, stochastic superresolution \citep{branicki2013dynamic}, and filtering Navier-Stokes equations with model error \citep{branicki2018accuracy}. It has been shown that the SPEKF model has much higher skill than classical Kalman filters using the so-called mean stochastic model (MSM) to capture the irregularity and intermittency in nature.

\begin{figure}[H]
\begin{center}
\includegraphics[scale=0.4]{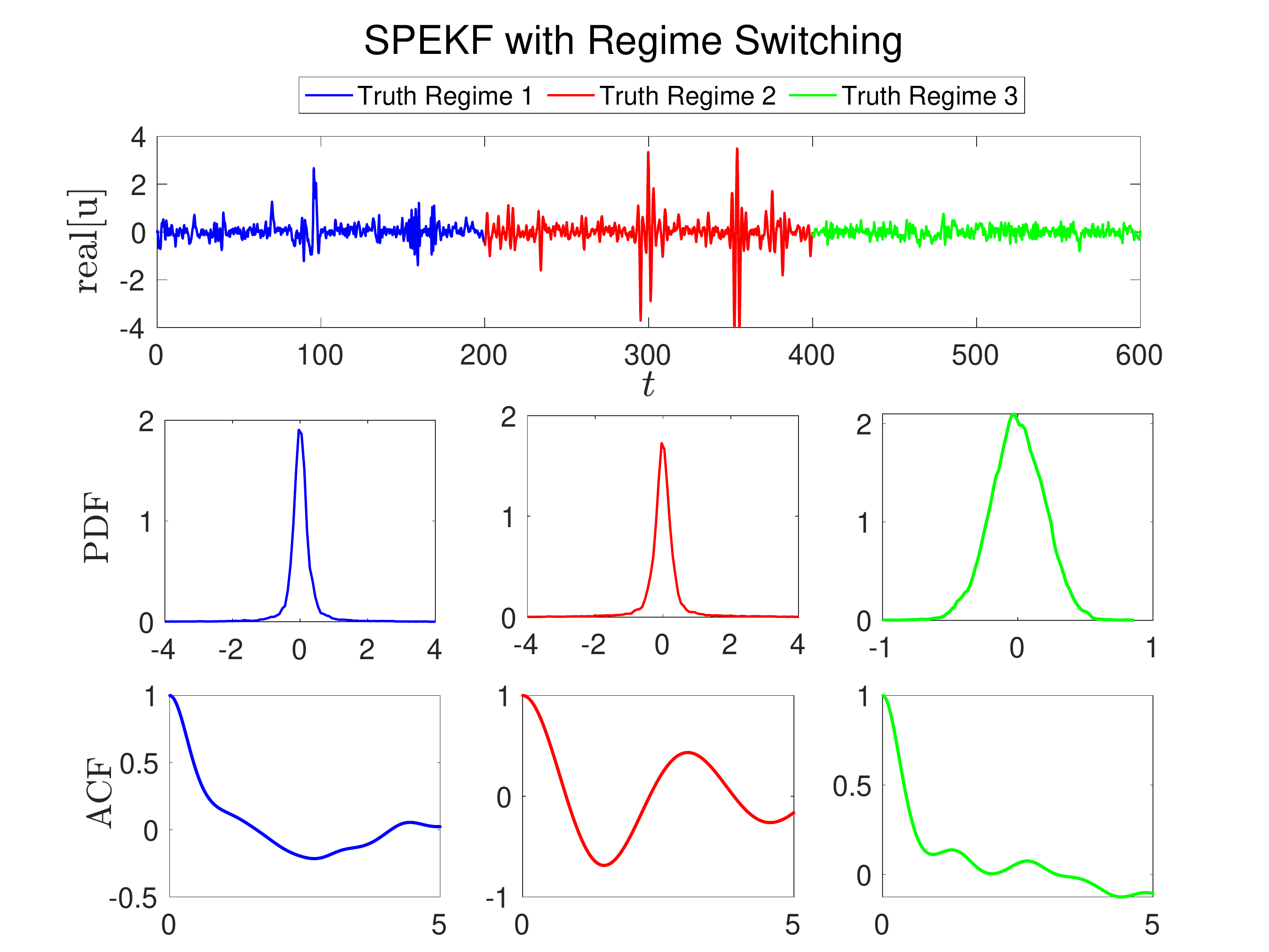}
\caption{The SPEKF model with three different regimes (twice the regime switching) as time evolves. In Regime
1, all three hidden variables $\gamma$, $\omega$, and $b$ have significant contributions to the dynamics of $u$. In Regime 2, $\omega$ is set to zero. In Regime 3, the stochastic damping $\gamma$ is removed
from the dynamics of u and the contribution from $\omega$ is added back. First row: the time series of the real part of the observed variable $u$. Second row: the PDF of $u$ in each regime. Third row: the ACF of $u$ in each regime.}
\label{fig: SPEKF_Property}
\end{center}
\end{figure}

In the following, the system will experience regime switching twice. In the starting regime (Regime 1), all three hidden variables have significant contributions to the dynamics of the observed process $u$. In Regime 2, $\omega$ will be set to zero; therefore, the stochastic phase does not influence the dynamics of the resolved variable $u$. Finally, in Regime 3, the stochastic damping $\gamma$ will be removed from the dynamics of $u$, but the contribution from $\omega$ will be added back. Figure~\ref{fig: SPEKF_Property} shows the observed trajectories and the associated statistics of the real part of the observed $u$ variable. Due to the stochastic damping, the time series in Regime 1 and Regime 2 are intermittent with multiple extreme events when the overall damping $\hat\gamma+\gamma(t)$ becomes positive. As a result, the PDFs are non-Gaussian fat-tailed. In contrast, very few extreme events are found in Regime 3, and the associated PDF is Gaussian. Similarly, the ACF in Regime 2 shows a clear oscillatory pattern when it decays. This is due to a dominant frequency of the time series coming from the constant phase $\hat\omega$. Such a regular oscillation in the ACF becomes less significant when the stochasticity is added to the phase term via $\omega(t)$.

The goal here is to detect the regime switching and reveal the dynamics in each regime. It is worth highlighting that the dimensions of the system in the three regimes are different. In the absence of $\omega$ or $\gamma$, the dimension reduces from $4$ to $3$ from Regime 1 to Regimes 2 and 3. However, such a change is unknown in practice and relies on the learning algorithm to detect it. In particular, there is no observed time series of $\gamma$, $\omega$, and $b$. Therefore, recovering these variables becomes an essential step in the online learning algorithm. To this end, the following procedure is adopted that incorporates data assimilation into the CEBoosting algorithm. Assume the model in Regime 1 is known. Each time when the new batch of time series of $u$ is obtained, such a model is utilized to sample a trajectory of the three unobserved variables $\gamma$, $\omega$, and $b$ conditioned on the observed trajectory of $u$. The sampled trajectory can be thought of as the analog of one ensemble member of the ensemble Kalman smoother solution of the system conditioned on the observed signal of $u$ \citep{evensen2000ensemble}, although the sampled trajectory using the SPEKF model can be written down using closed analytic formulae \citep{chen2022conditional}. See \citep{chen2020learning} for the implementation details of using such closed analytic formulae for data assimilation. This augments the unobserved components of each batch of data. Then the CEBooting algorithm is utilized to compute the causal relationships in light of the observed time series of $u$ and the sampled time series of $\gamma$, $\omega$, and $b$. If the causation entropy from any function involving $\gamma$ to the dynamics of $u$ is zero, then the $\gamma$ equation is eliminated from the final model structure. A similar logic applies to $\omega$ and $b$. Note that as $u$ and $b$ are complex-valued variables, the actual computation regards the real and imaginary parts as two processes in computing the causation entropy.

The main focus here is on identifying the dynamics of $u$. In particular, we aim to investigate if all three unobserved variables contribute to the observed process of $u$. To this end, the linear Gaussian models of $\gamma$, $\omega$, and $b$ are assumed to be fixed. A function library is designed with \{$u\gamma$, $u\omega$, $b$, $u^2$ , $\gamma\omega$ , $\gamma^2$, $\omega^2$, $b^2$, $\gamma b$, $\omega b$, $ub$\} to learn the structure of $u$.

\begin{figure}[H]
\centering
\includegraphics[scale=0.4]{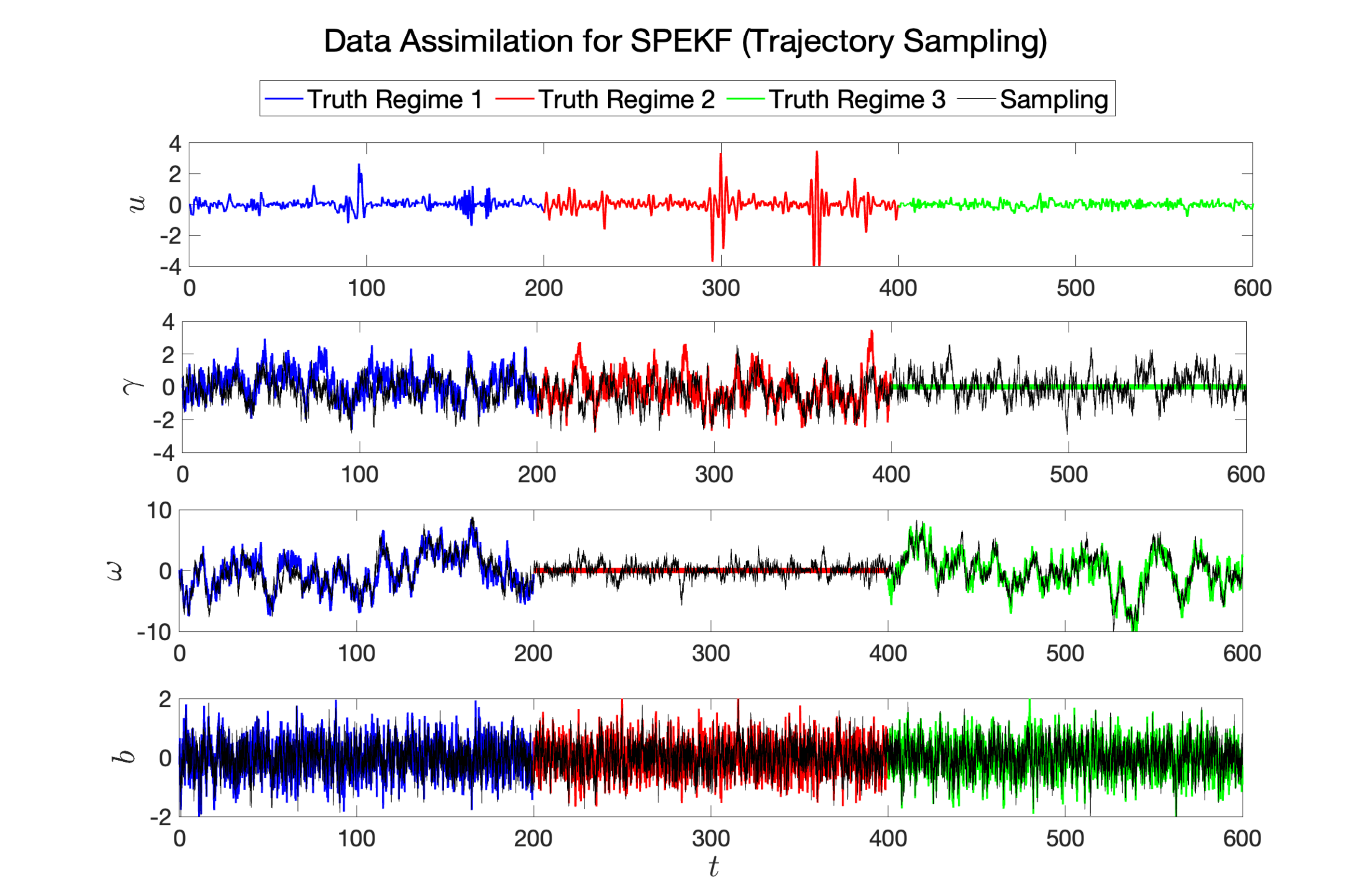}
\caption{Trajectory sampling for the SPEKF model with twice the regime switching. The black curve represents the conditional sampling of each hidden process given the observed variable $u$ in all three regimes.  }
\label{fig: SPEKF_DA_Sampling}
\end{figure}

Figure~\ref{fig: SPEKF_DA_Sampling} shows the real part of $u$ with regime switching and sampled trajectory of $\gamma$, $\omega$, $b$ across the three regimes. When the unobserved variables $\gamma$, $\omega$, or $b$ contribute to the dynamics of $u$, the sampled processes match the truth quite well. On the other hand, when $\omega$ and $\gamma$ disappear in Regimes 2 and 3, respectively, the sampling result provides random trajectories. The causation entropy from the terms involving these trajectories has no contribution to the dynamics of $u$. The first row of Table \ref{tab:SPEKF_CEM}, showing the causation entropies of Regime 1, is based on a time series with 200 time units. The second and the third rows show the causation entropies for the detected Regime 2 and Regime 3, respectively, using one batch of data with a length of 20 units. It is seen that the term $u\omega$ has a nearly zero causation entropy to $\dot{u}$ in Regime 2, where $\omega$ is sampled. Similarly, the sampled $\gamma$ leads to a nearly zero causation entropy from $u\gamma$ to $\dot{u}$ in Regime 3.

\begin{figure}[H]
\centering
\includegraphics[scale=0.4]{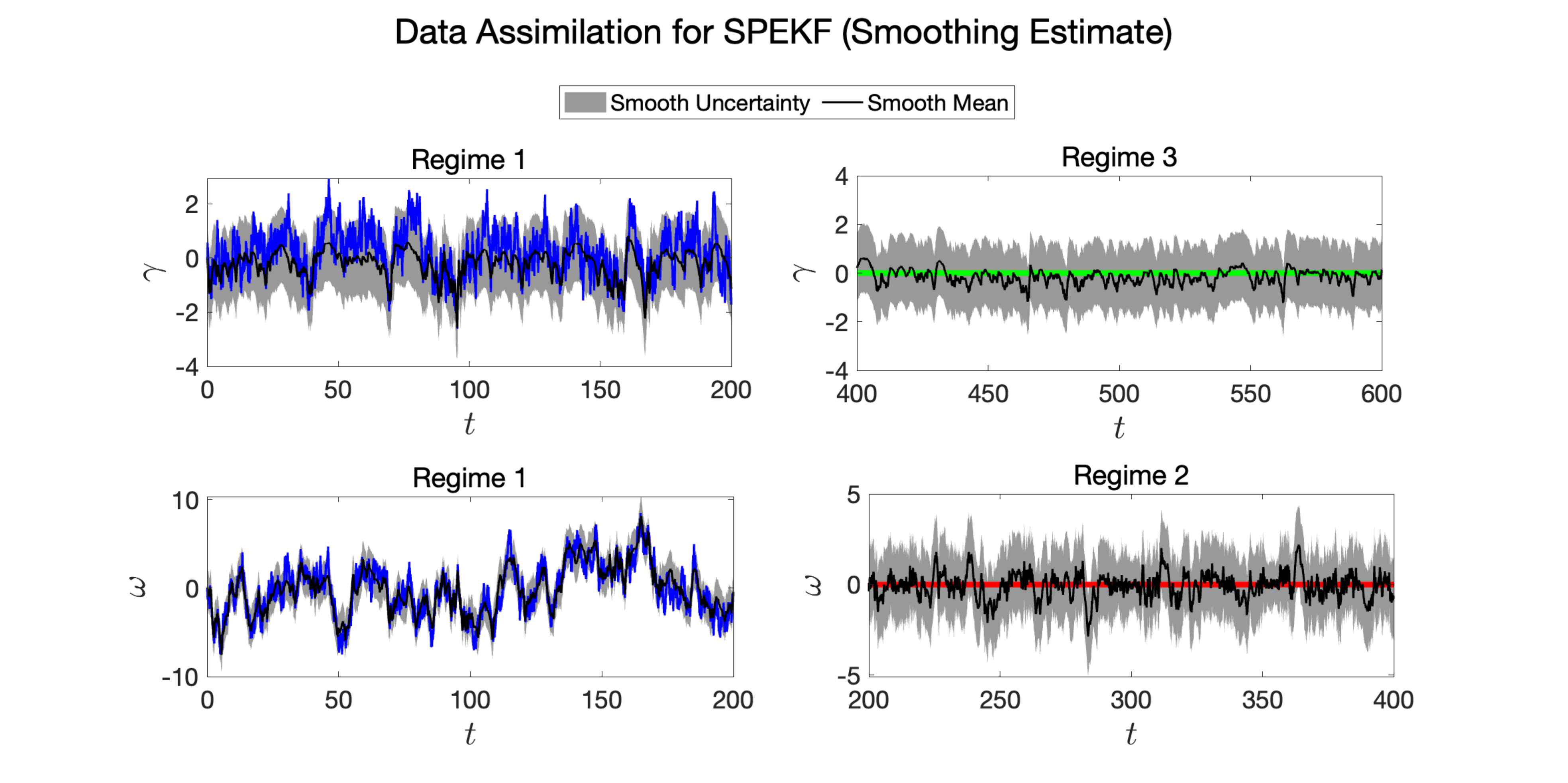}
\caption{The smoother mean and the smoother uncertainty (two standard deviations) of $\gamma$ and $\omega$ from the SPEKF model with twice the regime switching. The black curve shows the posterior mean of each hidden process conditioned on the observed $u$. The shaded gray area correspond to the 95\% confidence interval (two standard deviations) of the corresponding mean estimate.}
\label{fig: SPEKF_DA_Smoothing}
\end{figure}

Figure~\ref{fig: SPEKF_DA_Smoothing} shows the conditional mean and uncertainty (two standard deviations) of $\gamma$ and $\omega$ from the smoother solution. It is seen that when $\gamma$ and $\omega$ contribute to the dynamics of $u$ in Regime 1, the conditional mean follows the truth quite well. However, the conditional mean of $\omega$ and $\gamma$ behaves like a random trajectory in Regime 2 and Regime 3, respectively, with relatively large uncertainty. Thus, the sampled trajectories are formed from randomness and uncertainty, with no causal inference on the observed variable $u$.

\begin{table}[H]
\begin{center}
\caption{The SPEKF model - The CEM with Data Assimilation. First row: causation entropy with existing regime 1 data (200 time units). Second row: first batch (20 time units) causation entropy of regime 2. Third row: first batch (20 time units) causation entropy of regime 3.}
\label{tab:SPEKF_CEM}
\begin{adjustbox}{max width=1.\textwidth,center}
\begin{tabular}{ |c|c|c|c|c|c|c|c|c|c|c|c| }
 \hline
  & $u\gamma$ & $u\omega$ & $b$ & $u^2$ & $\gamma\omega$ & $\gamma^2$ & $\omega^2$ & $b^2$ & $\gamma b$ & $\omega b$ & $ub$ \\
 \hline
 $\dot{u}^{(1)}$ &\textbf{0.2527} &\textbf{0.5272}&\textbf{0.2289}&0.0077 & 0.0006&    0.0028&    0.0003&    0.0011&    0.0090&    0.0084&    0.0137 \\
 \hline
 $\dot{u}^{(2)}$     &\textbf{0.4525}&  0.0208 &    \textbf{0.1163}&    0.0491&    0.0070 &   0.0080 &   0.0019 &   0.0015&    0.0113 &   0.0035 &   0.0066
 \\
 \hline
 $\dot{u}^{(3)}$     & 0.0844 &   \textbf{0.7032} &    \textbf{0.3486}&    0.0048   & 0.0027 &   0.0019  &  0.0017&    0.0001&    0.0068  &  0.0106   & 0.0063
 \\
 \hline
 \end{tabular}
 \end{adjustbox}
\end{center}
\end{table}

\section{Conclusion}
\label{sec:conclusion}

Online nonlinear system identification with sequential data has recently become important in many applications, e.g., extreme weather events, climate change, and autonomous systems. In this work, we developed a causation entropy boosting (CEBoosting) framework for online nonlinear system identification. The CEBoosting algorithm aims to (i) discover a sparse residual model structure based on the aggregated causation entropy calculated from sequential data and (ii) calibrate the residual model with the identified sparse structure via least square estimation. If the true system experiences multiple regime switching, the proposed framework gradually identifies a summation of residual models, which has a close analogy to the statistical technique of boosting. We tested the proposed framework for complex systems with features including chaotic behavior, high dimensionality, intermittency and extreme events, and partial observations. The results show that the CEBoosting method can capture the regime switching and then calibrate residual models for various types of complex dynamical based on a limited amount of sequential data. It is worth noting that constraints can be naturally added to the learning algorithm. One important constraint is the so-called physics constraint \citep{majda2012physics}, which requires the total energy in the quadratic nonlinear terms to be conserved. It guarantees the long-term stability of the identified system. Such a constraint has yet to be explicitly incorporated into the current framework, although the resulting parameters in various non-Gaussian test cases shown in this work already roughly satisfy this constraint. Adding constraints can be easily achieved by imposing simple relationships between model parameters in the parameter estimation step, which still allows using closed analytic formulae for finding the parameters. See, for example, \citep{chen2020learning} for details. Other future work includes the uncertainty quantification of the aggregated causation entropy and the further study of other causality metrics.

\section*{Acknowledgments}
N.C. is partially funded by ONR N00014-19-1-2421 and ARO W911NF-23-1-0118. J.W. and C.C. are supported by the University of Wisconsin-Madison, Office of the Vice Chancellor for Research and Graduate Education with funding from the Wisconsin Alumni Research Foundation.

\clearpage
\appendix
\section{CEBoosting Algorithm}
\label{sec:algorithm}

Algorithm~\ref{alg:CEBoosting} presents the detailed procedures of the CEBoosting method, including (i) detecting regime switching, (ii) aggregating causation entropy matrix (CEM) until a consistent pattern of CEM is obtained, (iii) identifying a sparse model structure according to the aggregated CEM, and (iv) fitting the model parameters. The consistent pattern of CEM is determined by \eqref{eq:C_criterion}, i.e., the aggregated CEM does not change for $D$ data batches. In this work, we choose $D=4$ for all the numerical examples.

\begin{algorithm}[H]
\caption{Causation Entropy Boosting}
\begin{algorithmic}[1]
\State ${\fb}^* \gets \boldsymbol{\Xi}\boldsymbol{\Phi}$ \Comment{Current Model}
\For { $k=1,2,...$ }
    \State $\boldsymbol{\Phi} \gets \boldsymbol{\Phi}(\xb(t_{B_k}+m\Delta t))$
    \State $\dot{\xb} \gets [\xb(t_{B_k}+(m+1)\Delta t) - \xb(t_{B_k}+m\Delta t)]/\Delta t$
    \State $\mathbf{r} \gets \dot{\xb} - \boldsymbol{\Xi}\boldsymbol{\Phi}$ \Comment{Residual dynamics}
    \State $\mathbf{C} \gets C_{\phi_n \rightarrow r_i|[\Phi \setminus \phi_n]} \geq \overline{C}$
    \If{$\mathbf{C} =0$} \Comment{Same regime; Output current model}
        \State ${\fb}^* \gets \boldsymbol{\Xi}\boldsymbol{\Phi}$
    \Else \Comment{New regime detected}
        \State $K \gets 0$, $d \gets 1$, $D \gets 4$
        \While{ $d  \leq D-1 $} \Comment{Aggregate CEM until it is consistent for $D$ iterations}
            \State $K \gets K+1$
            \State $\mathbf{C} \gets \mathbf{C} + C^{(k+K)}_{\phi_n \rightarrow r_i|[\Phi \setminus \phi_n]}$
            \State $\mathbf{C}^+(K) \gets \frac{1}{K}\mathbf{C} \geq \overline{C}$ \Comment{Threshold aggregated CEM to get $\mathbf{C}^+$ with $0/1$ values}
            \If {$\mathbf{C}^+(K) = \mathbf{C}^+(K-1)$}
            \State $d \gets d+1$
            \Else
            \State $d \gets 1$
            \EndIf
        \EndWhile
        \State $K^* \gets K$ \Comment{Smallest $K^*$ that satisfies the criterion in \eqref{eq:C_criterion}}
        \State $\boldsymbol{\Xi}_r[\mathbf{C}^+(K^*) = 0] \gets 0$ \Comment{Select a sparse model structure and estimate parameters}
        \State $\boldsymbol{\Xi}_r[\mathbf{C}^+(K^*) = 1] \gets \argmin_{\boldsymbol{\Xi}_r} \sum_{m=1}^{MK^*}\|\mathbf{r}(t_{B_k}+m\Delta t) - \boldsymbol{\Xi}_r\boldsymbol{\Phi}(t_{B_k}+(m-1)\Delta t) \|^2$
        \State $\boldsymbol{\Xi} \gets \boldsymbol{\Xi} + \boldsymbol{\Xi}_r$
        \State ${\fb}^* \gets \boldsymbol{\Xi}\boldsymbol{\Phi}$
    \EndIf
    \State \Return ${\fb}^*$
\EndFor
\end{algorithmic}
\label{alg:CEBoosting}
\end{algorithm}

\end{document}